\documentclass[reqno]{amsart}

\usepackage{amsthm}
\usepackage{url}
\usepackage{amsmath}
\usepackage{amssymb}
\usepackage{amscd}
\usepackage{epsfig}
\usepackage{wrapfig}
\usepackage{graphicx}

\numberwithin{equation}{section}

\theoremstyle{plain}
\newtheorem{theorem}{Theorem}[section]
\newtheorem{proposition}{Proposition}[section]
\newtheorem{lemma}{Lemma}[section]
\newtheorem{corollary}{Corollary}[section]

\theoremstyle{definition}
\newtheorem{definition}{Definition}[section]
\newtheorem{counterexample}{Counterexample}[section]
\theoremstyle{remark}
\newtheorem{remark}{Remark}[section]
 
\newcommand{\sha}{\operatorname{Shadow}}

\newcommand{\E}{\mathbb E}
\newcommand{\Z}{\mathbb Z}
\newcommand{\R}{\mathbb R}
\newcommand{\Q}{\mathbb Q}
\newcommand{\N}{\mathbb N}
\newcommand{\dist}{\operatorname{dist}}
\newcommand{\const}{\operatorname{const}}
\newcommand{\Id}{\operatorname{Id}}
\newcommand{\Td}{\operatorname{Td}}
\renewcommand{\Pr}{\pi_{hg}}
\renewcommand{\S}{\mathbb S}
\begin{document}
\Large

\title{A.~D.~ALEXANDROV'S PROBLEM FOR BUSEMANN NON-POSITIVELY CURVED SPACES}

\author{P.~D.~Andreev}

\thanks{Supported by RFBR, grant 04-01-00315a}

\address{Pomor State University, Arkhangelsk, Russia}

\email{pdandreev@mail.ru}

\date{}
{\sloppy
\begin{abstract}{The paper is the last in the cycle devoted to the solution of Alexandrov's problem for non-positively curved spaces. Here we study non-positively curved spaces in the sense of Busemann. We prove that if $X$ is geodesically complete connected at infinity proper Busemann space, then it has the following characterization of isometries. For any bijection $f:X\to X$, if $f$ and $f^{-1}$ preserve the distance $1$, then $f$ is an isometry.\par
}
\end{abstract}

\maketitle

\markboth{{\upshape P.~D.~Andreev}}{{\upshape A.~D.~Alexandrov's problem}}

\vspace{3mm}

\tableofcontents

\section{Introduction}
The paper completes the cycle \cite{1}--\cite{3}, studying A. D. Alexandrov's problem for spaces with non-positive curvature. Previous papers were devoted to Alexandrov non-positively curved spaces. Here we deal with Busemann spaces defined in \cite{5}, see also \cite{6}--\cite{8}).

The main result of the paper is the following theorem.

\begin{theorem}\label{th1.1}
Let $(X, d_X)$ and $(Y, d_Y)$ be proper geodesically complete connected at infinity Busemann spaces, and $f: X \to Y$ be a bijection. Then the following statements are equivalent.
\begin{enumerate}
\item  The equality $d_X(x, y) = 1$ holds for points $x, y \in X$ iff $d_Y(f(x), f(y)) = 1$;
\item The inequality $d_X(x, y) \le 1$ holds for points $x, y \in X$ iff $d_Y(f(x), f(y)) \le 1$;
\item The inequality $d_X(x, y) < 1$ holds for points $x, y \in X$ iff $d_Y(f(x), f(y)) < 1$;
\item The map $f$ is an isometry of the space $(X, d_X)$ onto $(Y, d_Y)$.
\end{enumerate}
\end{theorem}

The trivial part of the theorem is the fact that statements (1)--(3) follows from (4). It is easy to observe that Theorem \ref{th1.1} has an equivalent formulation.

\begin{theorem}\label{th1.2}
Let the set $X$ is equipped with metrics $d_1$ and $d_2$, such that both spaces $(X, d_1), i = \overline{1, 2}$ satisfy the conditions of Theorem \ref{th1.1}. Then the following statements are equivalent.
\begin{enumerate}
\item The equality $d_1(x, y) = 1$ holds for points $x, y \in X$ iff $d_2(f(x), f(y)) = 1$; 
\item The inequality $d_1(x, y) \le 1$ holds for points $x, y \in X$ iff $d_2(f(x), f(y)) \le 1$; \item The inequality $d_1(x, y) < 1$  holds for points $x, y \in X$ iff  $d_2(f(x), f(y)) < 1$; 
\item The metrics $d_1$ and $d_2$ on the set $X$ coincide.
\end{enumerate}
\end{theorem}

Busemann curvature non-positivity condition is weaker than Alexandrov's one. Hence the class of Busemann spaces includes the subclass of $CAT(0)$-spaces. Some properties of $CAT(0)$-spaces are inherited in the considered class, the others undergone definite modifications in general.

The paper is organized as following. In Section \ref{1.05.08 12.27} we give formulations of necessary basic definitions and facts. The main part of the paper contains the proof of Theorem \ref{th1.2}. We use ideas and tools developed in earliest papers, adapted to the case of Busemann space. The proof is based on consideration of metrics $d_1$ and $d_2$ restricted to arbitrary straight line in the space $X$. There are several types of straight lines behavior. In Section \ref{1.05.08 12.28} we study the case of straight lines of higher rank and virtually higher rank, in Section \ref{1.05.08 12.30} the case of singular and virtually singular straight lines, and finally, in Section \ref{1.05.08 12.31} the case of strictly regular straight lines of strictly rank one. In the last section we present several counterexamples to the positive solution of A.D. Alexandrov's problem.

\section{Preliminaries}\label{1.05.08 12.27}

\subsection{Busemann non-positively curved spaces}
Let $(X,d)$ be a metric space. The ball with radius $\rho$ and center $x\in X$ is denoted  $B(x, \rho)$, the corresponding sphere $S(x, \rho)$. 

\begin{definition}
{\it A geodesic} in the space $(X,d)$ is by definition a locally homothetic map $c: I\rightarrow X$ where $I\subset \mathbb R$ is an interval or segment. The image of $I$ under the map $c$ is also called geodesic. The local homothety with a coefficient $\lambda > 0$ means that for any neighbourhood $U$ of arbitrary point  $t\in I$ the equality $d(c(s_1),c(s_2))= \lambda|s_1-s_2|$ holds for all $s_1,s_2 \in U$. The map $c$ presents {\it natural parameterization} of the geodesic if $\lambda = 1$ and {\it affine parameterization} or {\it parameterization proportional to natural} in general case. If $I = \mathbb R$, the geodesic $c$ is called {\it complete geodesic}. If the map $c$ is a homothety on the whole domain $I$, geodesic $c$ is called {\it minimizer}. In particular, minimizer defined on the segment $I=[a,b]\subset \mathbb R$ is called {\it segment} in the space $X$. It connects its {\it ends} $x = c(a)$ and  $y = c(b)$. The notation for the segment connecting points $x, y \in X$ is $[xy]$. {\it The straight line} is by definition a complete minimizer in $X$.  The space $(X,d)$ is called {\it geodesic} if any two its points can be connected by a segment. Geodesic space $X$ is called {\it geodesically complete}, if any geodesic in $X$ admits a continuation to a complete geodesic (not necessarily unique).
\end{definition}

\begin{definition}\label{14.03.08 9.30}
The geodesic space $X$ is called {\it Busemann non-positively curved} (or shortly {\it Busemann space}) if its metric is convex. This means the following. For any two segments $[xy]$ and $[x' y']$ with corresponding affine parameterizations $\gamma: [a, b] \to X$, $\gamma': [a', b'] \to X$, the function $D_{\gamma, \gamma'}: [a,b]\times [a',b'] \to \R$  defined by
$$D_{\gamma, \gamma'}(t, t') = |\gamma(t)\gamma'(t')|$$
is convex. Equivalently, the geodesic space $X$ is Busemann space if for any three points $x, y, z \in X$, the midpoint $m$ between $x$ and $y$ and the midpoint $n$ between $x$ and $z$ satisfy the inequality
\begin{equation}\label{busem}
|mn| \le \frac12 |yz|.
\end{equation}
Here the midpoint $m\in X$ between points $x, y$ satisfies equalities $d(x, m) = d(m, y) = \frac12 d(x, y)$.
\end{definition}

Busemann property of curvature non-positivity has a number of another equivalent formulations. The statements equivalent to Definition \ref{14.03.08 9.30} are listed in \cite[Proposition 8.1.2]{6}.  The simplest examples of Busemann spaces are $CAT(0)$-spaces and strictly convex normed spaces.

From now on the space $X$ satisfies to conditions of Theorem \ref{th1.2}. The distance between points $x, y \in X$ will be denoted $|xy|$.

\subsection{Normed strip Lemma}

Given a subset $A\subset X$ and $\epsilon>0$, the set
$$N_\epsilon(A):=\{x\in X\ |\ |xa| < \epsilon \text{\ for some\ } a\in A\}$$
is $\epsilon$-neighbourhood of $A$.

\begin{definition}
{\it Hausdorff distance} between closed subsets $A, B\subset X$ is by definition
$$d_H(A, B):=\inf \{\epsilon|A\subset N_\epsilon(B), B\subset N_\epsilon(A)\}.$$
In particular, if the value $\epsilon>0$ such that $A\subset N_\epsilon(B)$ and $B\subset N_\epsilon(A)$ does not exist, then $d_H(A, B) = \infty$.

Two straight lines $a, b : \mathbb R \to X$ are called {\it parallel}, if Hausdorff distance between them is finite:
$$d_H (a, b) < \infty.$$
\end{definition}

{\it Normed strip} ({\it Minkowski strip}) in the space $(X, d)$ is by definition a subset $L \subset X$ isometric to a strip between two straight lines in normed plane. The straight lines bounding the normed strip in $X$ are parallel.

The converse statement is also true. It is known as Rinow's normed strip lemma. We formulate the normed strip lemma as following.

\begin{lemma}[W. Rinow, \cite{Rin}, pp. 432, 463, \cite{Bow}, Lemma 1.1 and remarks]\label{lemmarinowa}
Every two parallel straight lines in Busemann space $X$ bound the normed strip in $X$.
\end{lemma}

\begin{remark} It is clear that Mincowski plane containing the strip isometric to normed strip in Busemann space is strictly convex.
\end{remark}

\begin{definition} We say that a straight line $a: \mathbb R \to X$ is {\it of higher rank} if it has parallel straight lines in $X$.
\end{definition}

\subsection{Compactifications of Busemann space}

The geometry at infinity of Busemann spaces have an essential difference from $CAT(0)$-spaces case. Two natural approaches to ideal compactification gives the same result in $CAT(0)$ case and can be different in the case of Busemann space. We refer for relations between two compactifications to \cite{10}. Here we only give necessary definitions and formulations.

\begin{definition}
The rays $c,d \colon [0, +\infty) \to X$ are called {\it asymptotic} if Hausdorff distance between them is finite:
$$\operatorname{Hd}(c,d) < + \infty.$$
The asymptoticity is an equivalence on the set of rays in $X$. The factor set $\partial_gX$ forms so called {\it geodesic ideal boundary} of $X$ and the union $\overline{X}_g = X \cup \partial_g X$ its {\it geodesic ideal compactification} of $X$. The topology on $\overline{X}_g$ called {\it cone topology} can be described as following. Given a basepoint $o \in X$ and a point $x \in \overline{X}_g$ we denote $[ox]$ a segment between them if $x \in X$, or a ray from $o$ to $x$ if $x \in \partial_gX$. By definition, the sequence $\{x_n\}_{n=1}^{+\infty}\subset \overline{X}_g$ converges to the point $x \in \overline{X}_g$ in the sense of the cone topology, if the sequence of natural patameterizations of segments (rays) $\{[ox_n]\}_{n=1}^{+\infty}$ converges to the natural parameterization of $[ox]$. In that case the convergence of parameterizations is uniform on common bounded subdomains of the parameter. Such defined cone topology on $\overline{X}_g$ does not depend on the choice of the basepoint $o$. Induced topology on the boundary $\partial_gX$ is also called {\it cone}. The identity map $\operatorname{Id}_X = i_g: X \to \overline{X}_g$ is an embedding of  $X$ to $\overline{X}_g$ as open dense subset $X = i_g(X) \subset \overline{X}_g$.
\end{definition}

We refer for the more complicated information about geometry of the boundary $\partial_gX$ to \cite{7} and \cite{8}.

Here we need some technical statement related to the cone topology on $\partial_gX$. Given a closed subset  $\mathcal{V} \subset \partial_g X$, point $o \in X$ and numbers $K, \varepsilon > 0$ we define \emph{$(o, K, \varepsilon)$-neighbourhood} of $\mathcal{V}$ as a set 
\[\mathcal{N}_{o, K, \varepsilon}(\mathcal{V}) := \{\zeta \in \partial_\infty X |\ \exists \xi \in \mathcal{V}, \ \zeta \in \mathcal{U}(\xi, o ,K, \varepsilon)\},\]
where
\[\mathcal{U}(\xi, o ,K, \varepsilon) := \{\eta \in \partial_\infty X |\  |\,c(K)d(K)|\,<\varepsilon, \  c=[o, \xi] ;\ d = [o, \eta] \}.\]

\begin{lemma}\label{12.04.08 8.04}
Fix a point $o \in X$ and a number $\varepsilon > 0$.
For any neighbourhood $\mathcal U$ of closed set $\mathcal V$ in the sense of the cone topology, there exists a number $K$ such that
$$V \subset \mathcal{N}_{o, K, \varepsilon}(\mathcal V) \subset \mathcal U.$$
\end{lemma}

\proof Assume to the contrary that for any $K > 0$ 
$$\mathcal{U}_{\xi_K, o, K, \varepsilon} \not\subset \mathcal U$$
for some $\xi_K \in V$. Fix a sequence $K_n \to +\infty$ and converging sequence $\{\xi_{K_n}\} \subset V$ of corresponding ideal points with neighbourhoods $\mathcal{U}_{\xi_{K_n}, o, K_n, \varepsilon} \not\subset \mathcal U$. We can do that because $V$ is compact. Denote 
$$\zeta = \lim\limits_{n \to \infty}\xi_{n}$$
and $K' > 0$ the number such that
$$\mathcal{U}_{\zeta, o, K', \varepsilon} \subset \mathcal U.$$
Then
$$\mathcal{U}_{\xi_{K_n}, o, 2K', \varepsilon} \subset \mathcal{U}_{\zeta, o, 2K', 2\varepsilon} \subset \mathcal U$$
for all but finitely many $n$. A contradiction to the choice of the sequence $\xi_{K_{n}}$. \qed

\begin{definition}
Let $(X, d)$ be a metric space $X$ and $C(X)$ be the space of continuous functions on $X$ with the topology of uniform convergence on bounded sets.  {\it Kuratowski embedding} $X \to C(X)$ is defined as following. Let $o \in X$ be a basepoint. Any point $x\in X$ is identified with {\it distance function} $d_x$ which acts by the formula
$$d_x(y) = |xy| - |ox|.$$
Let $C^\ast(X) = C(X)/\{\operatorname{consts}\}$ be a quotient space of $C(X)$ by the subspace of constants. Then the projection $p: C(X) \to C^\ast(X)$ generates embedding $\nu: X \to C^\ast(X)$ independent on the choice of the basepoint $o$. It is convenient to identify the space $X$ with its image $\nu(X)$. 

Let the space $X$ be proper and non-compact. {\it Horofunction compactification} of the space $X$ is by definition the closure of the image $\nu(X) \subset C^\ast(X)$. The horofunction compactification is denoted $\overline{X}_h$, the {\it  horofunction boundary} is $\partial_hX= \overline{X}_h \setminus X$. The map $\nu: X \to \overline{X}_h$ is embedding of $X$ to its horofunction compactification. Functions generating the horofunction boundary are called {\it horofunctions}. We think each horofunction as a limit of distance functions in the sense of the topology of uniform convergence on bounded sets. Given the horofunction $\Phi$, the corresponding point in the horofunction boundary is denoted $[\Phi] \in \partial_gX$. The important class of horofunctions in the Busemann space $X$ consists of {\it Busemann functions}. Every ray $c\colon \R_+ \to X$ generates corresponding Busemann function $\beta_c$ by the equality
$$\beta_c(y) = \lim_{t \to +\infty}(|yc(0)| - t).$$
Level sets of horofunctions are called {\it horospheres}, sublevels --- {\it horoballs}. The horosphere defined within the horofunction $\Phi$ by the equality $\Phi(x) = \Phi(x_0)$ where $x_0\in X$ is denoted $\mathcal{HS}(\Phi, x_0)$, the corresponding horoball is $\mathcal{HB}(\Phi, x_0)$.
\end{definition}

In \cite{11} M. Rieffel defined metric compactification of the space $X$. It is shown that metric compactification is equivalent to horofunction one.

The following theorem is proven in \cite{10}.

\begin{theorem}
Let $X$ be a proper non-compact Busemann space. Then there exists continuous surjection $\pi_{hg}: \overline{X}_h \to \overline{X}_g$ which coincide with the identification on $X$. If $\beta_c$ is Busemann function generated by the ray $c: \R_+ \to X$, then $\pi_{hg}([\beta_c])$ is a class of rays asymptotic to $c$ considered as a point in $\partial_gX$.
\end{theorem}

If $X$ is $CAT(0)$, the map $\pi_{hg}$ is a homeomorphism. From the other hand, the surjection $\pi_{hg}$ is not injective if $X$ is Minkowski space with singular norm. The preimage $\pi_{hg}^{-1}(\xi)$ consists of more than one point if $\xi$ corresponds to the singular direction of the norm.

\begin{definition} The point $\xi \in \partial_gX$ of geodesic ideal boundary is called  {\it regular} if its preimage $\pi_{hg}(\xi) \subset \partial_mX$ is one-point set. Otherwise the point  $\xi$ is {\it singular}. The straight line $a: \mathbb R \to X$ is called {\it regular} if both endpoints $a(-\infty)$ and $a(+\infty)$ are regular. Otherwise $a$ is called {\it singular}.
\end{definition}

It easily follows from the compactness of the space $\overline{X}_h$ and Hausdorffness of $\overline{X}_g$ that the map $\pi_{hg}$ is closed: the image of arbitrary closed subset in $\overline{X}_m$ is closed in $\overline{X}_g$. As a corollary, $\pi_{hg}$ satisfies to the following "weak openness" property.

\begin{lemma}\label{weakopenness}
For any point $\xi \in \partial_gX$ and any neighbourhood $\mathcal U$ of its preimage $\pi_{hg}^{-1}(\xi) \subset \partial_mX$ there exists a neighbourhood $\mathcal V$ of $\xi$ in $\partial_gX$ such that
$$\mathcal V \subset \pi_{hg}(\mathcal U).$$
\end{lemma}

\proof The image $\pi_{hg}(\partial_gX \setminus \mathcal U)$ is closed, so open subset
$$\mathcal V = \partial_gX \setminus \Pr(\partial_gX \setminus \mathcal U)$$
is demanded neighbourhood of $\xi$. \qed

\subsection{Virtual properties}

\begin{definition}
The finite collection of straight lines $a:=a_0,a_1,\dots, a_n:=b$ is called {\it asymptotic chain} if  for all $i = \overline{1,n}$ lines $a_{i-1}$ and $a_i$ are asymptotic in one of their directions. In that case we say that straight lines $a$ and $b$ are connected by the asymptotic chain. By definition, the straight line $b: \mathbb R \to X$ satisfies some property {\it virtually} if it is connected by asymptotic chain with the straight line $a$ which satisfies mentioned property. In further we need to consider virtually singular straight lines and straight lines virtually of higher rank. If the straight line $a$ is not a straight line virtually of higher rank (virtually singular), we say that $a$ is {\it strictly of rank one} ({\it strictly regular}).
\end{definition}

\subsection{Plan of the proof of Theorem \ref{th1.2}}

The equivalence of statements (1)--(3) in Theorem \eqref{th1.2} in the case of $CAT(0)$-space was proven in \cite{3}. Weakening the curvature conditions does not lead to changes in the proof. So we assume that the equivalence of statements (1)--(3) is proven. Our purpose is to show that these three claims imply the statement (4).

Consider a pair $x, y \in X$. By geodesic completeness of the space $X$ the segment $[xy]$ in the sense of the metric $d_1$ is contained in a straight line $a$ (not necessarily unique). We will prove that $a$ is a straight line in the sense of the metric $d_2$ as well, and metrics $d_1$ and $d_2$ are equal along $a$.

We need to study the following situations. The straight line $a$ can be of higher rank, virtually of higher rank or strictly of rank one. In the last case it can be singular, virtually singular or strictly regular. We prove the equality $d_1 = d_2$ along $a$ in all cases.

The main technique was developed in \cite{1}--\cite{3}. We use the notion of $r$-sequence introduced by V. Berestovski\v\i\ in \cite{1} and horospherical metric transfer from the straight line to its asymptotic straight line. Recall the definition of $r$-sequence following \cite{3}.

\begin{definition}
The homothety with coefficient $r > 0$\ $\mathbb Z \to X$ of integers $\Z$ to the space $X$ is called {\it $r$-sequence}.  We only consider the case $r = 1$ when the homothety becomes isometry, but we keep the term $r$-sequence for convenience. The segment of $r$-sequence $\{x_z\}_{z\in \mathbb Z}$ between $x_{z_{1}}$ and $x_{z_{2}}$ will be denoted
$$[x_{z_{1}}, x_{z_{1}+1}, \dots, x_{z_{2}}]_r.$$
Two $r$-sequences $\{x_z\}_{z \in \mathbb Z}$ and $\{y_z\}_{z \in \mathbb Z}$ are called {\it parallely equivalent} if Hausdorff distance between them is finite:
$$\operatorname{Hd}(\{x_z\}, \{y_z\}) < +\infty.$$
\end{definition}

The following result of V. Berestovski\v\i plays the crucial role in the proof in the case of $CAT(0)$-spaces (\cite{1}, Proposition 3.5). Let $X$ be $CAT(0)$-space that satisfies to conditions of Theorem \ref{th1.2}. Then the metric topology $\tau_m$ on $X$ is equal to the initial topology $\tau_f$ relative to the family of all  Busemann functions on $X$. We formulate corresponding proposition for the case of Busemann spaces as following.

\begin{proposition}
Let $X$ be geodesically complete connected at infinity proper Busemann space. Then the set of open horoballs corresponding to Busemann functions is a subbase for the metric topology on $X$.
\end{proposition}

\proof Given any ray $c = [x_0 \xi]$, the supplement $X \setminus \mathcal{HB}(\beta_c, x_0)$ of the closed horoball $\mathcal{HB}(\beta_c, x_0)$ is the union of open horoballs generated by Busemann functions. Indeed, 
$$X \setminus \mathcal{HB}(\beta_c, x_0) = \bigcup_{x \in \mathcal{HS}(\beta_c, x_0} \bigcup_{d \in r_x} \mathfrak{hb}(\beta_d, x).$$
Here $r_x$ denotes the set of rays $d: \R_+ \to X$ complement to the ray $[x\xi]$. The rest of the proof repeats the arguments of V. Berestovski\v\i\ from \cite{1}.\qed

The horospherical metric transfer is the procedure based on the following lemma. Its proof in \cite{3} does not change essentially after weakening the curvature conditions from $CAT(0)$ to Busemann spaces case.

\begin{lemma}\label{31.03.08 18.21}
Let the spaces $(X, d)$, $(X, d')$ satisfy conditions of theorem \ref{th1.2}.
Let images of maps $a, b:\R \to X$ be straight lines in $X$ with respect to both metrics $d$ and $d'$ and these straight lines are asymptotic in the direction of ideal point $\xi \in \partial_gX$ in the sense of metric $d$. If the equality $d = d'$ holds along $a$, then $d = d'$ along $b$ as well.
\end{lemma}

\subsection{The space of distances between asymptotic straight lines}

Fix an ideal point $\xi \in \partial_gX$ in the geodesic boundary of the space $X$. It defines the following pseudometric $\rho_\xi$ on $X$. For the points $x, y \in X$ put
$$\rho_\xi(x, y) = \operatorname{dist}([x\xi], [y,\xi]).$$
This means that
$$\rho_\xi(x, y) = \inf_{s, t \ge 0} |c(s), d(t)|,$$
where $c, d: \mathbb R_+ \to X$ are natural parameterizations of rays $[x_\xi]$ and $[y\xi]$ correspondingly.

The proof of the following claim is by direct checking of pseudometric axioms.

\begin{lemma}
The function $\rho_\xi$ is a pseudometric on $X$.
\end{lemma}

Denote $X_\xi$ the metric space obtained from $X$ within pseudometric $\rho_\xi$. The elements of $X_\xi$ are classes of points for which $\rho_\xi = 0$. We keep the notation $\rho_\xi$ for the metric on $X_\xi$. In particular, if the rays $c,d: \R_+ \to X$ are asymptotic in the direction $\xi$, the distance $\rho_\xi$ between their points is constant, and we denote this distance $\rho_\xi(c,d)$. In particular, the metric space $X_\xi$ may be one-point.

\begin{lemma}\label{2.04.08 21.10}
Let the rays $c, d : \mathbb R_+ \to X$ be asymptotic, $c(+\infty) = d(+\infty) = \xi$, $\beta_c$ and $\beta_d$ be corresponding Busemann functions. Then
$$0 \le \beta_c(d(0)) + \beta_d(c(0)) \le 2\rho_\xi(c, d).$$
\end{lemma}

\proof
Suppose that the ray $d' : \mathbb R_+ \to X$ in the direction of ideal point $\xi$ has common part with the ray $d$. We show the equality
\begin{equation}\label{18.03.08 20.45}
\beta_c(d(0)) + \beta_d(c(0)) = \beta_c(d'(0)) + \beta_{d'}(c(0)).
\end{equation}
In fact, if $d'(s) = d(t)$ for some $s, t \ge 0$, then
$$\beta_d(x) = \beta_{d'}(x) + t - s$$
for all  $x \in X$ and
$$\beta_c(d(0)) = \beta_c(d'(0)) - t + s.$$
Substitute $c(0)$ instead of $x$ in the first equality. Then adding of inequalities gives \eqref{18.03.08 20.45}.
In view of \eqref{18.03.08 20.45} we may assume that $\beta_c(d(0)) = 0$.

We claim that $\beta_d(c(0)) \ge 0$ in this case. Indeed, for any $\varepsilon > 0$ there exists  $T > 0$ such that for all $t> T$ 
$$|d(0)c(\tau)| < \tau + \varepsilon$$
and
$$\frac{t}{\tau} < \frac{\varepsilon}{|c(0)d(0)|}$$
for some   $\tau = \tau(t) > t$.
Let $p$ be the point of the segment $[d(0)c(\tau)]$ on distance
$$|d(0)p| = \frac{t}{\tau}\cdot |d(0)c(\tau)|$$
from  $d(0)$. We have
$$|pd(t)| \le \frac{t}{\tau}\cdot |c(\tau)d(\tau)| < \varepsilon $$
and $|pc(\tau)| < \tau - t + \varepsilon$. Consequently, from the triangle inequality
$$|c(0)p| > t - \varepsilon.$$
Hence
$$|c(0) d(t)| \ge |c(0) p| - |pd(t)| > t - 2 \varepsilon.$$
Taking into account arbitrariness of the choice of $\varepsilon$ and enlarging $t$ to infinity we obtain demanded estimation for $\beta_d(c(0))$ and for the sum $\beta_c(d(0)) + \beta_d(c(0))$ from below.
From the other hand, under supposing $\beta_c(d(0)) = 0$ take an arbitrary $\varepsilon > 0$ and  numbers $s, t > 0$ such that
$$|c(s)d(t)| < \rho_\xi(c,d) + \frac{\varepsilon}{4},$$
$$||c(s)d(0)| - s| < \frac{\varepsilon}{4}$$
and
$$||c(\Delta)d(t)| - t| < \frac{\varepsilon}{4}.$$
where $\Delta = \beta_d(c(0))$. Then triangle inequality gives
$$\Delta + t - \frac{\varepsilon}{4} < \Delta + |c(\Delta)d(t)| \le s + |c(s)d(t)| <  s +\rho_\xi(c,d) + \frac{\varepsilon}{4}$$
and
$$s - \frac{\varepsilon}{4} < |d(0)c(s)| \le t + |c(s)d(t)| < t + \rho_\xi(c,d) + \frac{\varepsilon}{4}.$$
Addition of the two inequalities gives
$$\Delta = \beta_c(d(0)) + \beta_d(c(0)) < 2 \rho_\xi(c,d) + \varepsilon.$$
Since $\varepsilon > 0$ was taken arbitrarily, we have necessary estimation from above for the sum $\beta_c(d(0)) + \beta_d(c(0))$. \qed 

Let $a : \mathbb R \to X$ be a straight lines with $a(+\infty) = \xi$ and $Y \subset X$ be a subset containing all points of straight lines parallel to $a$. Consider metric subspace $Y_\xi$ in the space $X_\xi$ obtained from $Y$.

\begin{lemma} \label{28.03.08 19.10}
The space $Y_\xi$ is Busemann non-positively curved space. It is one-point space iff $a$ is of rank one.
\end{lemma}

\proof
Let $b$ and $c$ be two straight lines parallel to $a$. Then $b$ parallel $c$ and they bound a normed strip $F \subset Y$. The strip $F$ is foliated by straight lines parallel to $a$ and it projects to a segment in the space $Y_\xi$. Consequently $Y_\xi$ is geodesic space. Let $b_\xi, c_\xi, d_\xi \in Y_\xi$ be three points obtained as projections to $Y_\xi$ of straight lines $b$, $c$ and $d$ correspondingly. Choose points $y\in c$ and $z\in d$ such that $|yz| = \rho_\xi(c_\xi, d_\xi)$. Also choose a point $x_1 \in b$ for which $|x_1 y| = \rho_\xi (b, c)$ and a point $x_2\in b$, for which $|x_2 z| = \rho_\xi(b, d)$.

\begin{figure}[ht]

\vspace{2cm}
\includegraphics[scale=0.75]{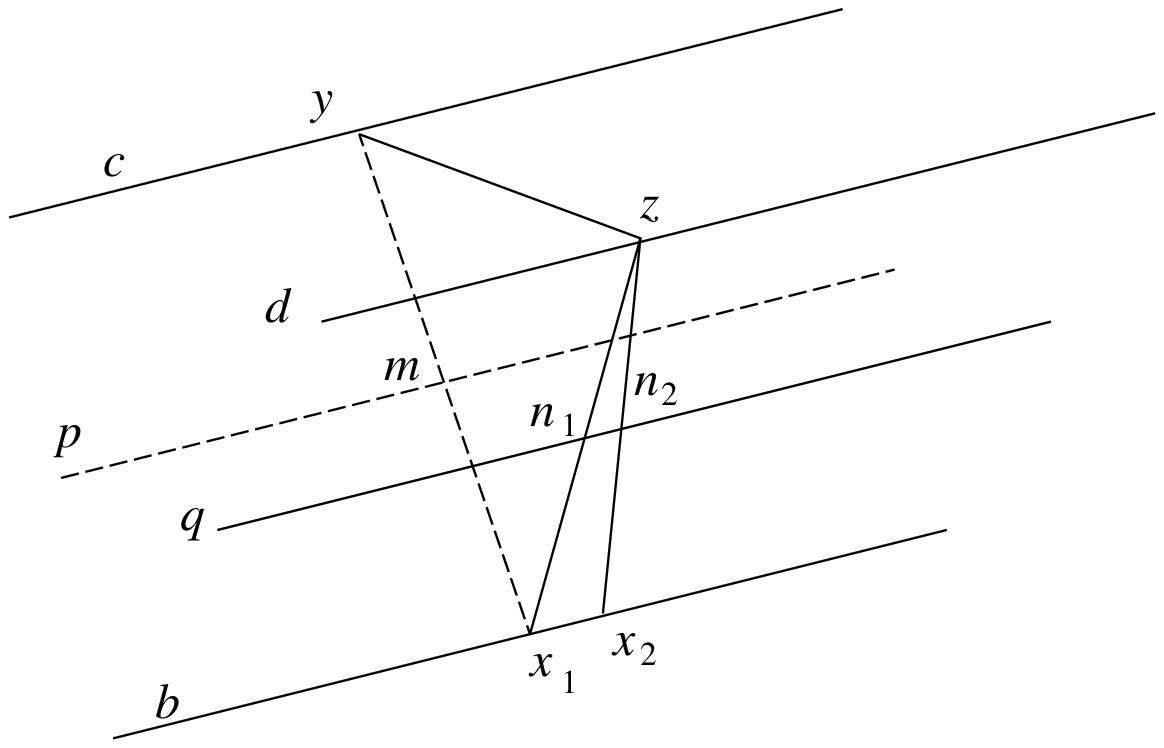}
\caption{}
\end{figure}
Let $m$ be the midpoint of the segment $[x_1y]$, $n_1$ the midpoint of the segment $[x_1z]$ and $n_2$ the midpoint of the segment $[x_2z]$, $p$ and $q$ be straight lines parallel to $a$ passing throw points $m$ and $n_2$ correspondingly, and $p_\xi$ and $q_\xi$ be their projections to $Y_\xi$. Then the straight line $q$ also passes throw the point $n_1$. Hence
$$\rho_\xi(p_\xi, q_\xi) \le |m n_1| \le \frac12|yz|.$$
Since $p_\xi$ are $q_\xi$ the midpoints  of the segments $[b_\xi c_\xi]$  and $[b_\xi d_\xi]$ in the space $Y_\xi$, the first claim of Lemma is proven. The second claim is obvious.

\section{Equality of metrics along the geodesic of higher rank}\label{1.05.08 12.28}

The main idea in consideration of straight lines of higher rank is inherited from the papers \cite{2} and \cite{3}. We study the construction of tapes introduced there. The only alteration is that we consider normed strips in the space $X$ instead of flat strips in $CAT(0)$ case. Such an alteration does not lead to essential changes in the proofs.

Recall the definition of $p$-tape.

\begin{definition}\label{28.03.08 19.06}
We say that the collection of $4p \quad (p \in \mathbb N)$ parallely equivalent $r$-sequences
\begin{equation}\label{29.03.08 20.32}
\{x_{i, j; z}\}_{z \in \mathbb Z},i= \overline{0, 3}, j= \overline{1, p}
\end{equation}  
forms a {\it $p$-tape}, if the following $4p+4$ points
$$
\begin{array}{l}
x_{i,\, 1,\, 0}, \dots, x_{i,\, p,\, 0}, \, i = \overline{0,\, 3} \\
x_{0,\, 1,\, 2p-1},\, x_{2,\, p,\, 1-2p},\, x_{3,\, p-1,\, 1-2p}, x_{3,\, p,\, 1-2p}
%x_{3,\, p-1,\, -p-3},\, x_{3,\, p-1,\, -p-2},\, x_{3,\, p-1,\, -p-1},
\end{array}
$$
generates in addition the system of segments of $r$-sequences:
$$
\left\{
	\begin{array}{l}
		[x_{0,\, 1,\, 0},\, x_{1,\, 1,\, 0},\, x_{2,\, 1,\, 0},\, x_{3,\, 1,\, 0}]_r \\
		%\left[x_{0,\, 2,\, 0},\, x_{1,\, 2,\, 0},\, x_{2,\, 2,\, 0},\, x_{3,\, 2,\, 0}\right]_r \\
		\dots \\
		\left[x_{0,\, p,\, 0},\, x_{1,\, p,\, 0},\, x_{2,\, p,\, 0},\, x_{3,\, p,\, 0}\right]_r \\
		%\left[x_{0,\, 1,\, 0},\, x_{1,\, p-1,\, -p-1},\, x_{2,\, p-1,\, -p-2},\, x_{3,\, p-1,\, -p-3}\right]_r \\
		\left[x_{0,\, 2,\, 0},\, x_{1,\, 1,\, 0},\, x_{2,\, p,\, 1-2p},\, x_{3,\, p-1,\, 1-2p}\right]_r\\
		\left[x_{0,\, 3,\, 0},\, x_{1,\, 2,\, 0},\, x_{2,\, 1,\, 0},\, x_{3,\, p,\, 1-2p}\right]_r\\
		\left[x_{0,\, 4,\, 0},\, x_{1,\, 3,\, 0},\, x_{2,\, 2,\, 0},\, x_{3,\, 1,\, 0}\right]_r\\
		\dots\\
		\left[x_{0,\, p,\, 0},\, x_{1,\, p-1,\, 0},\, x_{2,\, p-2,\, 0},\, x_{3,\, p-3,\, 0}\right]_r\\
		\left[x_{0,\, 1,\, 2p-1},\, x_{1,\, p,\, 0},\, x_{2,\, p-1,\, 0},\, x_{3,\, p-2,\, 0}\right]_r
	\end{array}
\right.\eqno{(2)}
$$
\end{definition}

\begin{figure}[h]
	\begin{picture}(360,60)
		\multiput(0,10)(80,0){3}{\line(5,1){120}}
		\multiput(40,34)(80,0){3}{\line(5,-1){120}}
		%\put(0,18){\line(5,-1){40}}
		\put(0,26){\line(5,-1){80}}
		%\put(0,18){\line(5,1){80}}
		\put(0,26){\line(5,1){40}}
		\put(240,10){\line(5,1){80}}
		%\put(360,10){\line(5,1){40}}
		%\put(320,34){\line(5,-1){80}}
		\put(280,34){\line(5,-1){40}}
		\multiput(0,10)(80,0){5}{\circle*{3}}
		\multiput(40,18)(80,0){4}{\circle*{3}}
		\multiput(0,26)(80,0){5}{\circle*{3}}
		\multiput(40,34)(80,0){4}{\circle*{3}}
		\put(0,0){$x_{0,\, 1,\, 0}$}
		\put(110,40){$x_{3,\, 1,\, 0}$}
		\put(225,0){$x_{0,\, p,\, 0}$}
		\put(300,0){$x_{0,\, 1,\, 2p-1}$}
		\put(18,40){$x_{3,\, p,\, 1-2p}$}
		\put(70,0){$x_{0,\, 2,\, 0}$}
		\put(70,35){$x_{2,\, 1,\, 0}$}
		\put(30,5){$x_{1,\, 1,\, 0}$}
		\put(260,40){$x_{3,\, p-1,\, 0}$}
		\put(265,5){$x_{1,\, p,\, 0}$}
	\end{picture}
\caption{$p$-tape.}
\label{tesma}
\end{figure}

We need the following technical statement. The notation $p - m - n - q$ means that points $m$ and $n$ belong to the segment $[pq]$ and divide this segment by three equal parts:
$$|pm| = |mn| = |nq| = \frac13|pq|.$$

\begin{lemma} \label{28.03.08 19.14}
Let $4p$ points $y_{ij}, i = \overline{0,3}, j = \overline{1,p}$ (not necessarily different) be given in Busemann space $Y$. Suppose that the following relations hold (see Fig. \ref{4.05.08 9.09}):
$$
\left\{
	\begin{array}{l}
		y_{01} - y_{11} - y_{21} - y_{31} \\
		\dots \\
		y_{0p} - y_{1p} - y_{2p} - y_{3p} \\
		y_{02} - y_{11} - y_{2p} - y_{3(p-1)} \\
		y_{03} - y_{12} - y_{21} - y_{3p} \\
		y_{04} - y_{13} - y_{22} - y_{31} \\
		\dots \\
		y_{0p} - y_{1(p-1)} - y_{2(p-2)} - y_{3(p-3)} \\
		y_{01} - y_{1p} - y_{2(p-1)} - y_{3(p-2)}
	\end{array}
\right.
$$
Then all points $y_{1j}$ coincide. The same is true for all points $y_{2j}$.
\end{lemma}

\begin{figure}[ht]

\vspace{2cm}
\hspace{-2cm}
\includegraphics[scale=0.8]{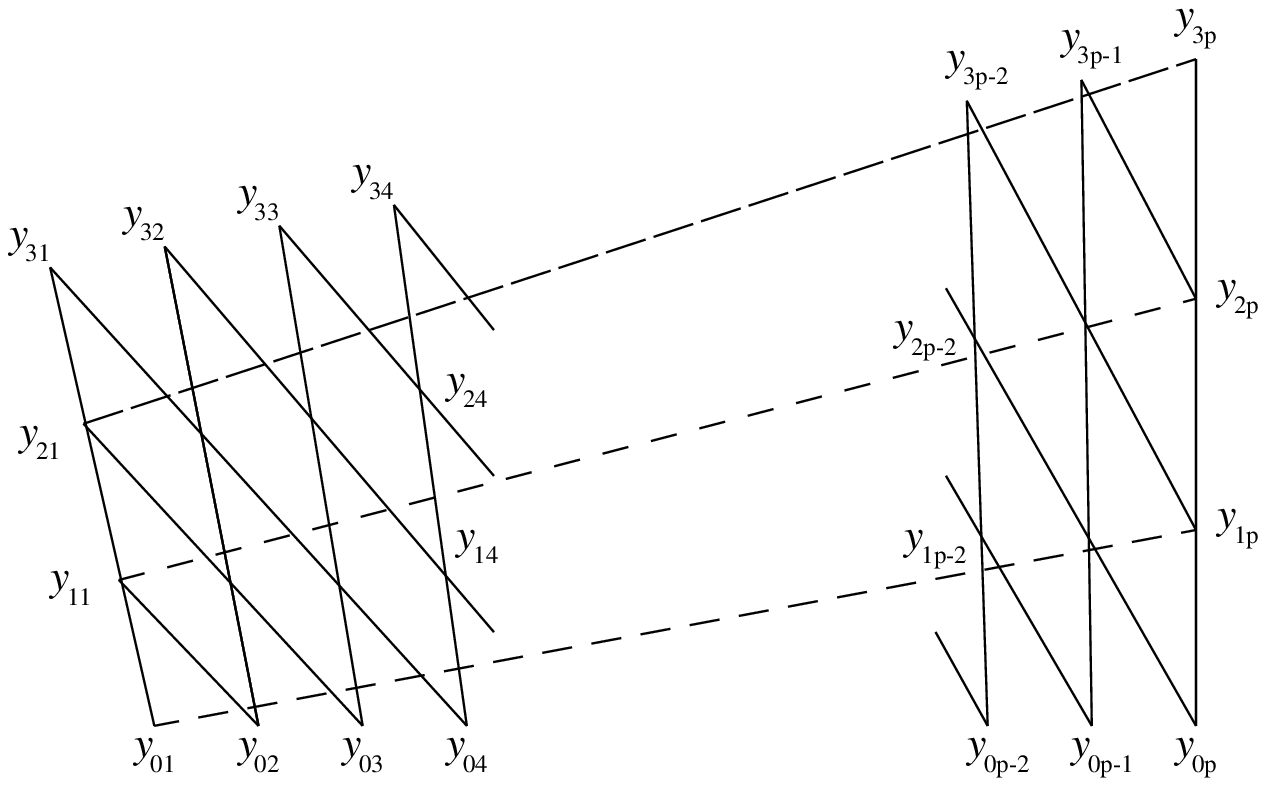}
\caption{}\label{4.05.08 9.09}
\end{figure}

\proof
Let $M$ be the maximum of distances $|y_{11}y_{12}|$, $|y_{12}y_{13}|, \dots$, $|y_{1(p-1)}y_{1p}|$, $|y_{1p}y_{11}|$, $|y_{21}y_{22}|$, $|y_{22}y_{23}|, \dots$, $|y_{2(p-1)}y_{2p}|$, $|y_{2p}y_{21}|$. After a renumeraton, if necessary, we may assume that $|y_{11}y_{12}| = M$. Then curvature non-positivity property with respect to the triangle $y_{02}y_{22}y_{2p}$ gives
$$|y_{2p}y_{22}| \ge 2M.$$ 
From the other hand, $$|y_{2p}y_{22}| \le |y_{2p}y_{21}| + |y_{21}y_{22}| \le 2M,$$
from where
$$|y_{2p}y_{22}| = 2M$$
and
$$|y_{2p}y_{21}| = |y_{21}y_{22}| = M.$$
Moreover, mentioned equalities mean that points $y_{2p}$, $y_{21}$ and $y_{22}$ belong to a straight line. Continuing in similar way, we get that all considering distances are $M$, all points $y_{1j}$ belong to a straight line and all points $y_{2j}$ also belong to another straight line. Such a configuration is possible only for $M = 0$. \qed

\begin{corollary}\label{30.03.08 13.00}
All $r$-sequences $x_{1, j, z}$ in Definition \ref{28.03.08 19.06} belong to one straight line. All $r$-sequences $x_{2, j, z}$ belong to one straight line as well.
\end{corollary}

\proof
Consider the set $Y$ formed by points of straight lines parallel to lines containing  $r$-sequences of the tape. By Lemma \ref{28.03.08 19.10} the space $Y_\xi$ is Busemann non-positively curved. Denote $y_{ij}$ the projection of $r$-sequence $\{x_{i, j, z}\}_{z \in \mathbb Z}$ to $Y_\xi$. Then points $y_{ij}$ form exactly the configuration described in Lemma \ref{28.03.08 19.14}. It follows that points $y_{1j}$ coincide. This point is the projection of one straight line containing $r$-sequences $\{x_{1, j, z}\}_{z \in \mathbb Z}$. Analogously, all $r$-sequences $\{x_{2, j, z}\}_{z \in \mathbb Z}$ lie in one straight line. \qed

Let $r$-sequences \eqref{29.03.08 20.32} form a $p$-tape in Busemann space $(X, d)$. Consider the segments $[x_{1, 1, 0}x_{3, 1, 0}]$ and $[x_{0, 2, 0}x_{2, 2, 0}]$. Their midpoints are $x_{2, 1, 0}$ and $x_{1, 2, 0}$ correspondingly and
$$|x_{2, 1, 0}x_{1, 2, 0}| = 1 = |x_{1, 1, 0}x_{0, 2, 0}| = |x_{3, 1, 0}x_{2, 2, 0}|.$$
Hence, given $t \in [0, 2]$, if $x_t \in [x_{1, 1, 0}x_{3, 1, 0}]$ is a point with $|x_{1, 1, 0}x_t| = t$ and $y_t  \in [x_{0, 2, 0}x_{2, 2, 0}]$ is a point with $|x_{0, 2, 0}y_t| = t$, then $|x_ty_t| = 1$. 

\begin{lemma}
The union $U$ of the segments $[x_ty_t]$ is convex subset in $X$ isometric to the parallelogram in the normed plane.
\end{lemma}

\proof For $t \in [0,2]$ and $s \in [0,1]$ denote $p_{(s,t)}$ the point of the segment $[x_ty_t]$ such that $|x_tp_{(s,t)}| = s$. Fix points $p_{(s_{1},t_{1})}$ and $p_{(s_{2},t_{2})}$; see Fig \ref{6.05.08 21.15}. For $\lambda \in [0, 1]$ denote $m_\lambda$ the point of the segment $[p_{(s_{1},t_{1})}p_{(s_{2},t_{2})}]$ such that 
$$|p_{(s_{1},t_{1})}m_\lambda| = \lambda |p_{(s_{1},t_{1})}p_{(s_{2},t_{2})}|.$$

\begin{figure}
	\begin{picture}(190,150)
		\put(20,50){\line(1,1){100}}
		\put(70,0){\line(1,1){100}}
		\put(50,80){\line(1,-1){50}}
		\put(60,90){\line(1,-1){50}}
		\put(100,130){\line(1,-1){50}}
		\put(60,70){\line(2,1){67}}
		\put(20,50){\circle*{3}}
		\put(70,0){\circle*{3}}
		\put(120,150){\circle*{3}}
		\put(170,100){\circle*{3}}
		\put(70,100){\circle*{3}}
		\put(120,50){\circle*{3}}
		\put(50,80){\circle*{3}}
		\put(100,30){\circle*{3}}
		\put(60,90){\circle*{3}}
		\put(110,40){\circle*{3}}
		\put(100,130){\circle*{3}}
		\put(150,80){\circle*{3}}
		\put(60,70){\circle*{3}}
		\put(127,103){\circle*{3}}
		\put(73,77){\circle*{3}}
		\put(-10,50){$x_{1,1,0}$}
		\put(40,0){$x_{0, 2, 0}$}
		\put(45,105){$x_{2,1,0}$}
		\put(125,45){$x_{1,2,0}$}
		\put(95,155){$x_{3,1,0}$}
		\put(175,95){$x_{0, 2, 0}$}
		\put(38,59){$p_{(s_{1},t_{1})}$}
		\put(130,110){$p_{(s_{2},t_{2})}$}
		\put(32,79){$x_{t_{1}}$}
		\put(100,20){$y_{t_{1}}$}
		\put(80,73){$m_\lambda$}
		\put(83,132){$x_{t_{2}}$}
		\put(150,70){$y_{t_{2}}$}
		\put(0,92){$x_{(1-\lambda)t_{1}+\lambda t_{2}}$}
		\put(110,30){$y_{(1-\lambda)t_{1}+\lambda t_{2}}$}
		\end{picture}
\caption{}\label{6.05.08 21.15}
\end{figure}

It follows from the convexity of the metric $d$ that
$$1 = |x_{(1-\lambda)t_1 + \lambda t_2}y_{(1-\lambda)t_1 + \lambda t_2}| \le |x_{(1-\lambda)t_1 + \lambda t_2} m_\lambda| + |m_\lambda y_{(1-\lambda)t_1 + \lambda t_2}| \le
$$
$$\le \left((1-\lambda)|x_{t_{1}}p_{(s_{1},t_{1})}| + \lambda |x_{t_{2}}p_{(s_{2},t_{2})}|\right) +
\left((1-\lambda)|p_{(s_{1},t_{1})}y_{t_{1}}| + \lambda |p_{(s_{2},t_{2})}y_{t_{2}}|\right) \le$$
$$\le (1-\lambda)|x_{t_{1}}y_{t_{1}}| + \lambda |x_{t_{2}}y_{t_{2}}| = 1.$$
Since the left and right sides of the inequality above coincide, all the inequalities must be equalities. It follows that $m_\lambda$ is $p_{((1-\lambda)s_1 + \lambda s_2, (1-\lambda)t_1 + \lambda t_2)}$. Hence the subset $U$ is convex. Moreover, all the maps $\lambda \to p_{((1-\lambda)s_1 + \lambda s_2, (1-\lambda)t_1 + \lambda t_2)}$ represent the affine parameterizations of the segments $[p_{(s_{1},t_{1})}p_{(s_{2},t_{2})}]$.

\begin{figure}[b]
	\begin{picture}(190,150)
		\put(70,0){\line(1,1){120}}
		\put(70,0){\circle*{3}}
		\put(110,40){\circle*{3}}
		\put(150,80){\circle*{3}}
		\put(190,120){\circle*{3}}
		\put(110,40){\line(-2,1){80}}
		\put(150,80){\line(-2,1){80}}
		\put(30,80){\line(4,1){160}}
		\put(70,0){\line(0,1){120}}
		\put(70,60){\circle*{3}}
		\put(70,120){\circle*{3}}
		\put(30,80){\circle*{3}}
		\put(110,100){\circle*{3}}
		\put(78,0){$p_{(2s_{0}-s,t_{0})}$}
		\put(115,38){$p_{(s_{0},t_{0})}$}
		\put(155,78){$p_{(s,t_{0})}$}
		\put(195,118){$p_{(2s-s_{0},t_{0})}$}
		\put(-5,88){$p_{(s_{0} + \delta,t_{0} + \sigma)}$}
		\put(40,130){$p_{(s + \delta,t_{0}+ \sigma)}$}
		\end{picture}
\caption{}\label{6.05.08 22.08}
\end{figure}

Given arbitrary $\delta\in (-2, 2)$ and $\sigma \in (-1 ,1)$, consider the function $\rho_{\delta, \sigma}$ defined in the appropriate part of the rectangle $[0,1] \times [0,2]$ by the equality
$$\rho_{\delta, \sigma}(s, t) = d(p_{(s,t)}, p_{(s+\delta, t+\sigma)}).$$
We claim that the function $\rho_{\delta, \sigma}$ is constant on its domain of representation $$\mathcal D_{\delta, \sigma} \subset [0, 2] \times [0, 1]$$
and
\begin{equation}\label{7.05.08 8.17}
\rho_{\lambda\delta, \lambda\sigma}(s, t) = \lambda \rho_{\delta, \sigma}(s,t)
\end{equation}
for all $\lambda \in [0,1]$. It is sufficient to prove the claim for a small neighbourhood of the arbitrary interior point $(s_0, t_0) \in \mathcal D_{\delta, \sigma} \cap (0, 2) \times (0,1)$.

We prove that if  $(2s - s_0, t_0), (2s_0 - s, t_0) \in \mathcal D_{\delta, \sigma}$, then $\rho_{\delta, \sigma}(s, t_0) = \rho_{\delta, \sigma}(s_0, t_0)$. Indeed, the Busemann inequality for the triangle $p_{(2s - s_{0},t_{0})}p_{(s_{0},t_{0})}p_{(s_{0}+ \delta, t_{0}+ \sigma)}$ gives (see Fig. \ref{6.05.08 22.08})
$$|p_{(s,t)}p_{(s+\frac12\delta,t_0+\frac12\sigma)}| \le \frac12 |p_{(s_{0},t_{0})}p_{(s_{0}+ \delta,t_{0}+\sigma)}|$$
and consequently
$$\rho_{\delta\sigma}(s,t_0) \le \rho_{\delta\sigma}(s_0, t_0).$$
From the other hand, the Busemann inequality for the triangle $p_{(2s_{0}-s, t_{0})}p_{(s,t)}p_{(s+\delta,t_0 +\sigma)}$ gives
$$|p_{(s_{0},t_{0})}p_{(s_{0}+\frac12\delta, t_{0}+ \frac12 \sigma)}| \le \frac12 |p_{(s,t0}p_{(s+\delta,t_0+\sigma)}|$$
and
$$\rho_{\delta\sigma}(s_0, t_0) \le \rho_{\delta\sigma}(s, t_0).$$
Hence the equality holds
$$\rho_{\delta\sigma}(s_0, t_0) = \rho_{\delta\sigma}(s, t_0).$$
Similarly we obtain
$$\rho_{\delta\sigma}(s_0, t_0) = \rho_{\delta\sigma}(s_0, t)$$
for all $t \in [0, 1]$ such that $(s_0, 2t-t_0), (s_0, t-2t_0) \in \mathcal D_{\delta, \sigma}$. As a corollary, the equality holds
$$\rho_{\delta\sigma}(s_1, t_1) = \rho_{\delta\sigma}(s_2, t_2)$$
for all pairs $(s_1, t_1), (s_2, t_2) \in \mathcal D_{\delta, \sigma}$. The equality \eqref{7.05.08 8.17} is obvious.

Now, define the norm $N$ in the affine plane $A^2$ with coordinates $(\alpha, \beta)$ by the equality 
$$N(\alpha, \beta) = \frac{1}{\lambda}\rho_{|\lambda \alpha|, |\lambda \beta|}(p_{st}),$$
where $\lambda > 0$, $|\lambda \alpha| \le 2$, $|\lambda \beta| \le 1$ and $(s, t) \in \mathcal D_{|\lambda \alpha|, |\lambda \beta|}$ are taken arbitrarily. The previous consideration shows that the norm $N$ does not depend on the choice of $\lambda$ and $(s, t)$ satisfying mentioned conditions. It follows from the convexity of the metric in $X$ that the normed space $(A^2, N)$ is strictly convex. By the definition of the norm $N$, the parallelogram $[0, 2] \times [0, 1] \subset A^2$ is isometric to $U$. \qed

\begin{lemma}\label{30.03.08 14.32}
Let $r$-sequences \eqref{29.03.08 20.32} form a $p$-tape and the straight line $a: \R \to X$ contains points $x_{1, j, z}$ for $1 \le j \le p$ and $z \in \Z$. Then 
\begin{equation}\label{30.03.08 14.33}
x_{1, j, z} = a\left(\frac{(j-1)(2p-1)}{p} + z\right)
\end{equation}
for all $j \in \{1, \dots, p\}$ and $z \in \Z$.
\end{lemma}

\proof

It follows from the previous lemma that
$$|x_{1, 1, 0}x_{1, 2, 0}| = |x_{2, 1, 0}x_{2, 2, 0}|.$$ 
Analogously, the equalities hold
$$|x_{1, j-1, 0}x_{1, j, 0}| = |x_{2, j-1, 0}x_{2, j, 0}|$$
for all $j \in \{2, \dots, p\}$,
$$|x_{2, j, 0}x_{2, j+1, 0}| = |x_{1, j+1, 0}x_{1, j+2, 0}|$$
for all $j \in \{1, \dots, p-2\}$ and
$$|x_{2, p-1, 0}x_{2, p, 0}| = |x_{1, p, 0}x_{1, 1, 2p-1}|.$$
Since from the Corollary \ref{30.03.08 13.00} the points $x_{1, j, 0}$  belong to one straight line containing also $x_{1, 1, 2p-1}$, these points divide the segment $[x_{1, 1, 0}x_{1, 1, 2p-1}]$ to $p$ equal parts. The rest of the proof is obvious. \qed

Now suppose that the map $a: \R \to X$ represents a straight line of higher rank in the metric space $(X, d_1)$. Then the image $a(\R)$ can lie in the interior of some normed strip or in the other case $a(\R)$ is boundary line of any normed strip containing it. Suppose the first option. Let $\mathcal F$ be a normed strip containing $a(\R)$ in its interior.

\begin{lemma} \label{30.03.08 11.04}
There exists a number $P > 0$ such that for all natural $p > P$ the normed strip $\mathcal F$ contains a $p$-tape as in Definition \ref{28.03.08 19.06} with 
$$x_{1, j, z} = a\left(\frac{(j-1)(2p-1)}{p} + z\right).$$
\end{lemma}

\proof Let a number $L > 0$ be such that the strip $\mathcal F$ contains a substrip of width $3L$ with boundary lines on the distances $L$ and $2L$ from $a$. Fix a point $q \in \mathcal F$ with $|a(0)q| = 1$ and
$$0 < \dist(q, a) < \min\{1, L\}.$$
then there exists number $t$ with $0 < |t| < 2$ such that $|a(t)q| = 1$. Since the metric of $\mathcal F$ is strictly convex, the number $t$ and the point $a(t)$ are defined uniquely. Take a number $P > 0$ such that $2/P < 2-|t|$. It is easy to see that $P$ satisfies the claim. \qed

\begin{lemma}\label{30.03.08 11.26}
Let the collection of $r$-sequences \eqref{29.03.08 20.32} forms $p$-tape in the sense of metric $d_1$. Then it also forms $p$-tape in the sense of the metric $d_2$.
\end{lemma}

\proof Follows immediately from conditions on metrics $d_1$ and $d_2$ and Definition \ref{28.03.08 19.06}. \qed

Now we are ready to prove the equality of metrics along the straight line $a$ of higher rank in the case when a passes in the interior of some normed strip $\mathcal F$.

\begin{lemma} \label{31.03.18.06}
Let the map $a: \R \to X$ represent a straight line the sense of the metric $d_1$, and $a$ passes in the interior of the normed strip $\mathcal F$. Then the map $a$ represents a straight line in the sense of the metric $d_2$ and we have equality $d_1 = d_2$ along it.
\end{lemma}

\proof By Lemma \ref{30.03.08 11.04} there exists a number $P > 0$ such that for all natural $p > P$ the normed strip $\mathcal F$ contains $p$-tape defined by the collection \eqref{29.03.08 20.32} with
$$a\left(\frac{k}{p}\right) = x_{1, j, z},$$
where
$$j = 1 - p \cdot \left\{\frac{k}{p}\right\},$$
$\left\{\frac{k}{p}\right\}$ denotes fractional part of $k/p$, and
$$ z = \frac{k - (j-1)(2p-1)}{p}.$$
By Lemma \ref{30.03.08 11.26}, every such $p$-tape is also $p$-tape in the sense of the metric $d_2$. By Corollary \ref{30.03.08 13.00} points $x_{1, j, z}$ belong to the image of the map $a': \R \to X$ representing natural parameterization of a straight line in the sense of the metric $d_2$, and by Lemma \ref{30.03.08 14.32} they satisfy equalities \eqref{30.03.08 14.33}. Taking different natural values $p$, we obtain that $a'(q) = a(q)$ for all rational $q$. Since metrics $d_1$ and $d_2$ are equivalent, the equality $a'(t) = a(t)$ holds for all $t \in \R$. Hence the claim. \qed

\begin{corollary}\label{31.03.08 18.19}
Let the map $a: \R \to X$ be the natural parameterization of a straight line of higher rank in the sense of the metric $d_1$. Then the map $a$ is also natural parameterization of a straight line in the sense of the metric $d_2$.
\end{corollary}

\proof If the image $a(\R)$ passes in the interior of a normed strip, the result is proven in Lemma \ref{31.03.18.06}. If $a(\R)$ is boundary straight line of normed strip $\mathcal F$, it is the limit of naturally parameterized straight lines. Since the metrics $d_1$ and $d_2$ are equivalent, $a(\R)$ is the limit of naturally parameterized straight lines in the sense of the metric $d_2$ as well. Hence the claim. \qed

Finally we have the result.

\begin{theorem}\label{17.04.08 0.25}
Let $a$ be a straight line virtually of higher rank in the sense of metric $d_1$. Then $a$ is straight line virtually of higher rank in the sense of the metric $d_2$ and metrics coincide along it:
$$d_1(x, y) = d_2(x, y)$$
for all $x, y \in a$.
\end{theorem}

\proof Follows immediately from Corollary \ref {31.03.08 18.19} and Lemma \ref{31.03.08 18.21}. \qed

\section{Equality of metrics along singular straight line}\label{1.05.08 12.30}

\subsection{Double spherical transfer}

In this section, we prove that metrics $d_1$ and $d_2$ coincide along singular straight line of rank one. The singularity of the straight line $a: \mathbb R\to X$ means that at least one of enpoints $a(+\infty)$ or $a(-\infty)$ is singular point of the geodesic ideal boundary $\partial_gX$. Since $a$ is of rank one in the sense of metric $d_1$, it follows that the image $a(\R)$ is also straight line of rank one in the sense of the metric $d_2$ (see \cite{3} for details). Hence we only need to prove the equality $d_1 = d_2$. We assume that the singular ideal point is $\xi = a(+\infty) \in \partial_gX$.

We need to discuss some properties of horofunctions and Busemann functions now.

\begin{lemma}
Let $c = [o\xi]$ be the ray and $\beta_c$ be corresponding Busemann function. Let $\Phi$ be a horofunction with $\Phi(o) = 0$ and $\Pr([\Phi]) = \xi$. Then $\beta_c \ge \Phi$: 
$$\beta_c(x) \ge \Phi(x)$$
for all $x \in X$.
\end{lemma}

\proof It is sufficient to show that if $\Phi(y) = 0$, then $\beta_c(y) \ge 0$. Fix a point $y \in X$ with $\Phi(y) = 0$ and arbitrarily small number $\varepsilon > 0$. Represent the horofunction $\Phi$ as a limit function for the sequence of distance functions $d_{x_{n}}$. Here
$x_n \to \xi$ in the sense of the cone topology on $\partial_gX$. Let the number $K$ be such that $$|\beta_c(y) - (|yc(t)| - t)| < \frac{\varepsilon}{4}$$
for all $t \ge K$ and the number $N$ be such that 
$$|d_{x_{n}}(y)| < \frac{\varepsilon}{2}$$
and
$$|c(K)\sigma_n(K)| < \frac{\varepsilon}{4}$$
for all $n \ge N$. Here $\sigma_n: [0, |ox_n|] \to X$ denotes the natural parameterization of the segment $[ox_n]$. Now from the triangle inequality we obtain
$$|yx_n| \le |yc(K)| + |c(K)\sigma_n(K)| + |\sigma_n(K)x_n| < \left(\beta_c(y) + K + \frac{\varepsilon}{2}\right) + (|ox_n| - K) <$$
$$< \beta_c(y) + |ox_n| + \frac{\varepsilon}{2}.$$
Hence
$$0 = \Phi(y) < |yx_n| - |ox_n| + \frac{\varepsilon}{2} < \beta_c(y) + \varepsilon.$$
Since $\varepsilon$ is taken arbitrarily small, it follows the claim. \qed

\begin{lemma}\label{2.04.08 20.43}
Let $\xi \in \partial_g X$ be a singular point of the boundary $\partial_gX$. Then the set $\Pr^{-1}(\xi) \subset \partial_m X$ contains more than one Busemann function.
\end{lemma}

\proof Consider Busemann function $\beta_c$ generated by the ray $c = [o\xi]$ and the horofunction $\Phi \ne \beta_c$ with $\Phi(o) = 0$ and $\Pr([\Phi]) = \xi$.  Take a point $y$ where $\Phi(y) = 0 < \beta_c(y)$. Consider the ray $d=[y\xi]$ and corresponding Busemann function $\beta_d$. We have $\beta_d(y) = \Phi(y) = 0$. Consequently
$$\beta_d(x) \ge \Phi(x) = 0,$$
$$\beta_c(x) - \beta_d(x) \le 0$$
and
$$\beta_c(y) - \beta_d(y) \ge 0.$$
Hence the difference $\beta_c - \beta_d$ is not constant and points $[\beta_c], [\beta_d] \in \partial_hX$ are different. \qed

\begin{corollary}\label{2.04.08 22.37}
Given a ray $a: \mathbb R_+ \to X$ with $a(+\infty) = \xi \in \partial_g X$, where $\xi$ is a singular point of geodesic ideal boundary $\partial_gX$, there exists a ray $b \colon [0, +\infty) \to X$ asymptotic to $a$, such that the difference of Busemann functions $\beta_a - \beta_b$ is non-constant. Moreover, the ray $b$ can be chosen so that $\beta_a(a(0)) = \beta_a (b(0))$ and $\beta_b(a(0)) \ne \beta_b(b(0))$.
\end{corollary}

\proof The first claim follows immediately from Lemma \ref{2.04.08 20.43} and the second claim from the first one. \qed

\begin{definition}
Let $a, b: \R \to X$ be asymptotic straight lines with common endpoint at infinity $\xi = a(+\infty) = b(+\infty) \in \partial_gX$. Let $\beta_a$ (correspondingly, $\beta_b$) be Busemann function defined from the ray $a|_{\R_+}$ (correspondingly $b|_{\R_+}$). Double horospherical transfer $T_{a \leftrightarrow b}: a \to a$ is defined by the condition: $T_{a \leftrightarrow b}(x) = x' \in a$ if $\beta_b(x') = \beta_b(y)$, where $y \in b$ is a point such that $\beta_a(y) = \beta_a(x)$. In other words, if $x = a(t)$, then $x' = a(t')$, where $t' - t = \beta_a(b(0)) + \beta_b(a(0))$. It follows from Lemma \ref{2.04.08 21.10} that $t - t' \ge 0$.
\end{definition}

\begin{remark}
It is clear that $T_{a \leftrightarrow b}$ is isometric translation of straight line $a$. If $T_{a \leftrightarrow b}(x) = x$ for some point $x \in a$, then $T_{a \leftrightarrow b} = \Id(a)$. In particular, this holds when the point $\xi \in \partial_gX$ is regular.
\end{remark}

\begin{theorem}
Let $a: \mathbb R \to X$ be a singular straight line in metric space $(X, d_1)$. Then the map $a$ represents singular straight line in metric space $(X, d_2)$ as well, and $d_1(x, y) = d_2(x, y)$ for any pair of points $x, y \in a$.
\end{theorem}

\proof
The case of the straight line virtually of higher rank was studied in previous section. So we may think the straight line $a$ to be strictly of rank one. 
Also we may think the ideal point $\xi = a(+\infty)$ to be singular.
Since the straight line $a$ is singular, it admits asymptotic straight line $b$ with
$b(+\infty) = a(+\infty) = \xi$ 
such that the difference of Busemann functions
$\beta_a $ and $\beta_b$ defined by rays
$a|_{\mathbb R_+}$ and $b|_{\mathbb R_+}$ correspondingly is non-constant. Since horospheres corresponding to Busemann functions in 
the sense of metric $d_1$ are also horospheres in the sense of metric $d_2$, it follows that the straight line $a$ is singular in the metric $d_2$ as well. 
By the corollary \ref{2.04.08 22.37} we can choose naturally parameterized line $b$ so that
$\beta_a(b(0)) = 0$ and $\beta_b(a(0)) > 0$.
Consider a segment
$[a(0)b(0)]$. 
When the point $x$ moves by this segment continuously from $b(0)$ to $a(0)$, the function
$\beta_a(x)$ is non-positive by the convexity of the horoball
$\mathcal{HB}(\beta_a, a(0))$. 
Define the following function
$\mathcal{B}(x)$. 
Let
$c_x = [x \xi] : \mathbb R \to X$
be a ray from
$c_x(0) = x$ in the direction of the point
$c_x(+\infty) = \xi$. 
Denote
$\beta_{c_{x}}$ corresponding Busemann function. 
The function $\mathcal{B}(x)$ is defined by the equality
$$\mathcal{B}(x) = \beta_{c_{x}}(a(0)).$$
By Lemma \ref{2.04.08 21.10} the value $\mathcal{B}$ depends continuously on the point of the segment $[a(0)b(0)]$. Hence the sum $\beta_a(x) + \mathcal{B}(x)$ is continuous when $x$ moves in the segment from $b(0)$ to $a(0)$ and it takes values from $\beta_b(a(0)) > 0$ to $0$. In particular, there exists a natural number $N \in \mathbb N$, such that for all natural $n > N$ 
$$\beta_a(x_n) + \mathcal{B}(x_n) = \frac1n$$
for some point $x_n \in [a(0)b(0)]$.
We denote arbitrary straight line containing the ray $c_x$ by the same symbol $c_x$. The double horospherical transfer  $T_{a \leftrightarrow c_{x_{n}}}$ maps the point $a(0)$ to $a(1/n)$. So
\begin{equation}\label{7.04.08 21.49}
(T_{a \leftrightarrow c_{x_{n}}})^n(a(0)) = a(1)
\end{equation}
and the equality \eqref{7.04.08 21.49} holds in the sense of both metrics $d_1$ and $d_2$. Consequently,
$$d_1(a(0), a(t)) = d_2(a(0), a(t))$$
for any rational $t \in \mathbb Q$.
Moreover, 
\begin{equation}\label{7.04.08 21.54}
d_1(a(t_1), a(t_2)) = d_2(a(t_1), a(t_2))
\end{equation}
for any $t_1, t_2 \in\R$ with $t_2 - t_1 \in \Q$.
Since metrics $d_1$ and $d_2$ are topologically equivalent and have common the incidence relation on straight lines of higher rank, we conclude that the equality \eqref{7.04.08 21.54} is true for any values $t_1, t_2 \in \mathbb R$. \qed

\begin{corollary}\label{17.04.08 0.27}
Let $a$ be virtually singular straight line in metric space $(X, d_1)$. Then it is virtually singular in metric space $(X, d_2)$ and $d_1(x, y) = d_2(x, y)$ for any $x, y \in a$.
\end{corollary}

\section{Equality of metrics on strictly regular straight line strictly of rank one}
\label{1.05.08 12.31}
\subsection{Tits relations on the boundary $\partial_gX$}
The main tool for the proof of equality for metrics $d_1$ and $d_2$ along the straight line $a$ in the case when $a$ is strictly regular and strictly of rank one is scissors defined in \cite{3}. The principle of the proof also does not change essentially. But in addition, we need to study metric properties of the boundary at infinity in the case of Busemann space. In particular, we can not use Tits metric $\Td$ on $\partial_gX$ because this metric admits no general definition with properties of Tits metric in the case of $CAT(0)$-space. 

In \cite{12}, we introduced a collection of binary relations that can be considered as substitute of Tits distance. There are two key values of Tits distance: $\pi$ and $\pi/2$. The most of geometric applications of Tits metric is based on the comparison of Tits distance between ideal points with these key values. But inequalities of type $\Td(\xi, \eta) > \pi$ etc. have purely geometric description without using Tits metric itself. This allows to introduce the following trick. We define the collection of binary relations on ideal boundaries $\partial_gX$ and $\partial_hX$ corresponding to comparison of Tits distance with $\pi$ and $\pi/2$. Here we only need relations of type $\Td(\xi, \eta) > \pi$ and $\Td(\xi, \eta) \le \pi$. Recall the definition (\cite{12}, Definition 3.2).

\begin{definition}\label{deftits}
Let $(X, o)$ be a pointed proper Busemann space.
Let rays $c, d: \R_+ \to X$ with common beginning
$$c(0) = d(0) = o$$ 
represent points $\xi  = c(+\infty)$ and $\eta  = d(+\infty)$ in the boundary $\partial_gX$. 
The function $\delta_o: \partial_gX \times \partial_gX \to [0, \pi]$ is well-defined by the equality
$$\delta_o(\xi, \eta)  =  \lim_{t \to +\infty}\delta_{o, \xi, \eta}(t),$$
where
$$\delta_{o, \xi, \eta}(t)  = \frac{|c(t)d(t)|}{2t}.$$

Given ideal points $\xi, \eta \in \partial_gX$ we define the following binary relations:
\begin{itemize}
\item $\Td(\xi, \eta) < \pi$ if $\delta_o(\xi, \eta) < \pi$
\item $\Td(\xi, \eta) \le \pi$, if for any neighbourhoods $U(\xi)$ and $V(\eta)$ of this points in the sense of cone topology on $\partial_gX$ there exist points $\xi'  \in U(\xi)$ and $\eta'  \in V(\eta)$ with $\Td(\xi' , \eta' ) < \pi$;
\item $\Td(\xi, \eta) \ge \pi$, if $\Td(\xi, \eta) < \pi$ does not hold;
\item $\Td(\xi, \eta) > \pi$ if $\Td(\xi, \eta) \le \pi$ does not hold;
\item $\Td(\xi, \eta)  = \pi$ if $\Td(\xi, \eta) \ge \pi$ and $\Td(\xi, \eta) \le \pi$ hold simultaneously.
\end{itemize}
\end{definition}

One of the main consequence of the definition above is the following theorem.

\begin{theorem}[\cite{12}, Theorem 3.1]\label{Td>pi}
Let $X$ be a proper Busemann space.
If $\Td(\xi, \eta) > \pi$, then there exists a geodesic $a : \R \to X$ with ends $a(-\infty) = \xi$ and $a(+\infty) = \eta$.
\end{theorem}

\begin{corollary} \label{7.05.08 21.21} Given a straight line $a: \R \to X$ of rank one,  endpoints $\xi = a(+\infty)$ and $\eta = a(-\infty)$ have cone neighbourhoods $\mathcal U_+ = \mathcal U(\xi)$ and $\mathcal U_- = \mathcal U(\eta)$ such that for any pair of ideal points $\zeta \in \mathcal U_+$ and $\theta \in \mathcal U_-$ 
$$\Td(\zeta, \theta) > \pi$$
and there exists a straight line $b: \R \to X$ with endpoints $b(+\infty) = \zeta$ and $b(-\infty) = \theta$.
\end{corollary}

Another application of Definition \ref{deftits} is the following criterion for the existence of normed half planes with given boundary.

\begin{definition}
Normed half plane in the space $X$ is by definition the subspace isometric to a half plane in Minkowski plane.
\end{definition}

\begin{theorem}[\cite{12}, Theorem 3.2]
Let $X$ be a proper Busemann space. Given a geodesic $a: \R \to X$ with endpoints $\xi  = a(+\infty)$ and $\eta  = a(-\infty)$ passing throw $a(0)  = o$, the following conditions are equivalent.
\begin{enumerate}
\item $\Td(\xi, \eta) = \pi$; 
\item there exist horofunctions $\Phi$ centered in $\xi$ and $\Psi$ centered in $\eta$, for which the intersection of horoballs
\begin{equation}\label{intersechoroballs}
\mathcal{HB}(\Phi, o) \cap \mathcal{HB}(\Psi, o)
\end{equation}
is unbounded;
\item  $a$ bounds a normed half plane in $X$.
\end{enumerate}
\end{theorem}

In the connection with the cone topology we need the following property of strictly regular straight lines strictly of rank one.

\begin{lemma}\label{15.04.08 19.31}
Let $a: \R \to X$ be strictly regular straight line strictly of rank one with endpoints $\xi = a(+\infty)$ and $\eta = a(-\infty)$. Then for any $\varepsilon > 0$ there exists cone neighbourhoods $\mathcal U_+$ of $\xi$ and $\mathcal U_-$ of $\eta$ such that the following holds. If a straight line $b: \R \to X$ has endpoints $b(+\infty) \in \mathcal U_+$ and $b(-\infty) \in \mathcal U_-$, then
$$|a(0) b(t)| < \varepsilon$$
for some $t \in \R$.
\end{lemma}

\proof 
Since $a$ has rank one, then $\Td(\xi, \eta) > \pi$ and the ideal points $\xi$, $\eta$ have cone neighbourhoods $\mathcal{U}'_+$ and $\mathcal{U}'_-$ correspondingly, such that if $\zeta \in \mathcal{U}'_-$ and $\theta \in \mathcal{U}'_+$, then $\Td(\zeta, \theta) > \pi$. Consequently, the points $\zeta$ and $\theta$ admits the straight line $c: \R \to X$ with $c(-\infty) = \zeta$ and $c(+\infty) = \theta$. In that case the straight line $c$ can be connected with $a$ by the asymptotic chain $a, b, c$ where $b(-\infty) = a(-\infty) = \eta$ and $b(+\infty) = c(+\infty) = \theta$.

Denote $\beta_+$ and $\beta_-$ Busemann functions defined from rays $[a(0)\xi]$ and $a(0)\eta]$ correspondingly. Then, since $a$ is of rank one,
$$\beta_+(x) + \beta_-(x) \ge 0$$
and equality holds if and only if $x \in a$. By convexity of functions $\beta_{\pm}$, there exist numbers $\delta_1 > 0$ such that
$$\beta_+(x) + \beta_-(x) >  \delta_1$$
for all $x \in X \setminus B(a(0), \varepsilon/2)$, and $\delta_2 > \delta_1$ such that
$$\beta_+(x) + \beta_-(x) >  \delta_2$$
for all $x \in X \setminus B(a(0), \varepsilon)$.
Denote 
$$\mu = \min\left\{\frac{\delta_1}{2}, \frac{\delta_2-\delta_1}{2}\right\}$$
and
$$\mathcal{U}_{\pm}^\ast = \{g \in C(X, \R) \, |\, \forall x \in B(a(0), \varepsilon) |g(x) - \beta_{\pm}(x) - \delta_1/2| < \mu$$
the neighbourhoods of the functions $\beta_{\pm}+ \delta/2$ in the space $C(X, \R)$.
If $g_+ \in \mathcal U_+^\ast$ and $g_- \in \mathcal U_-^\ast$, then
$$g_+(a(0) + g_-(a(0)) < 0$$
and
$$g_+(x) + g_-(x) > 0$$
for all $x \in S(a(0), \varepsilon)$. The projections $p|_{\mathcal{U}_{\pm}^\ast}: \mathcal{U}_{\pm}^\ast \to C^\ast(X)$ contain neighbourhoods in $\partial_hX$ of points $[\beta_+], [\beta_-] \in \partial_hX$. If $[\Phi] \in p(\mathcal{U}_+^\ast$, then some horofunction $\Phi \in [\Phi]$ belongs to $\mathcal{U}_+^\ast$. Hence the intersection of horoballs $\mathcal{HB}(\Phi, a(0)) \cap \mathcal{HB}(\beta_-, a(0))$ is containing in $B(a(0), \varepsilon)$ and compact. There exists a straight line $b: \R \to X$ with $b(-\infty) = \eta$ and $b(+\infty) = \pi_{hg}([\Phi])$. Since $a$ is strictly regular, the point $\pi_{hg}([\Phi])$ is regular and $\Phi$ is Busemann function.  Then the projection $\pi_{hg}|_{p|_{\mathcal{U}_{+}^\ast}}$ is a homeomorphism to some neighbourhood $\mathcal{U}''_+$ of the point $\xi \in \partial_gX$. Analogously, projection $\Pr: p|_{\mathcal{U}_{-}^\ast} \to \partial_gX$ is homeomorhism onto some neighbourhood $\mathcal{U}''_{-}$.

Denote $\mathcal{U}_{\pm} = \mathcal{U}'_{\pm} \cap \mathcal{U}''_{\pm}$. By the construction, neighbourhoods $\mathcal U_+$ and $\mathcal U_-$ satisfy the claim of the Lemma. \qed

\subsection{Scissors}

\begin{definition}[\cite{3}, Definition 4.1]
We say that straight lines $a, b, c, d: \R \to X$ in the space $X$ form {\it scissors} with {\it center} $x\in X$ if
\begin{itemize}
\item $a(-\infty)=b(-\infty)$; 
\item $a(+\infty)=c(+\infty)$; 
\item $c(-\infty)=d(-\infty)$;  
\item $b(+\infty)=d(+\infty)$; и
\item $b \cap c = x$.
\end{itemize}

We denote the configuration of scissors $\langle a,b,c,d ;x \rangle$ (fig. \ref{scissors}). The straight line $a$ is called \emph{base} of scissors. In the case when the straight line $a$ is strictly regular, the four ideal points serving as endpoints of straight lines $a, b, c, d$ generates exactly four classes of horofunctions, namely four classes of Busemann functions presented by functions $\beta_{a(\pm \infty)}$ and $\beta_{d(\pm \infty)}$ with
$$\beta_{a(\pm \infty)}(x)=\beta_{d(\pm \infty)}(x) = 0.$$
\end{definition}

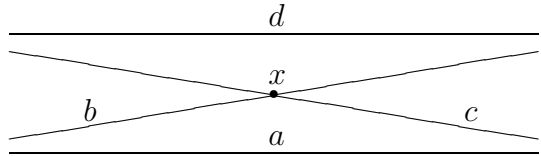
\begin{figure}[h]
\begin{picture}(200,50)
\put(0,0){\line(1,0){200}}
\put(0,5){\line(6,1){200}}
\put(200,5){\line(-6,1){200}}
\put(0,45){\line(1,0){200}}
\put(100,22){\circle*{3}}
\put(98,26){$x$}
\put(98,3){$a$}
\put(98,48){$d$}
\put(172,12){$c$}
\put(28,12){$b$}
\end{picture}
\caption{Scissors $\langle a,b,c,d ;x \rangle$}
\label{scissors}
\end{figure}

When scissors are fixed, they generates a translation $T$ of the base $a$ as following. Let  $R_{ac}$ be horospherical transfer from the straight line $a$ to $c$ generated by Busemann function $\beta_{a(+\infty)}$: the point $m \in a$ moves to the unique point $m^\prime = R_{ac}(m) \in c$ for which $\beta_{a(+\infty)}(m^\prime) = \beta_{a(+\infty)}(m)$. Analogously, one defines transfers $R_{cd}$, $R_{db}$ and $R_{ba}$, which are isometric maps of corresponding straight lines. Note that the transfers above are defined independently on the choice of the metric $d_1$ or $d_2$ on the space $X$ and they act isometrically with respect to both metrics.

\begin{definition}
\emph{The translation} $T$ is a composition
\[T:=R_{ba}\circ R_{db} \circ R_{cd} \circ R_{ac} : a \to a.\]
Obviously, $T$ is an isometry of the straight line $a$ preserving its direction. \emph{The shift} $\delta T$ of the translation $T$ defined as the difference 
$$\delta_T = \beta_{a(-\infty)}(T(m)) - \beta_{a(-\infty)}(m)$$ 
is independent on the choice of the point $m\in a$.
\end{definition}

The quantity $\delta T$ can be described as following. Let $\beta_{a-}$, $\beta_{a+}$, $\beta_{d-}$ and $\beta_{d+}$ be Busemann functions with centers $a(\pm \infty)$ and $d(\pm \infty)$ correspondingly, such that there exist points $p \in a$ and $q \in d$, for which $\beta_{a-}(p) = \beta_{a+}(p) = 0$ and $\beta_{d-}(q) = \beta_{d+}(q) = 0$. 

\begin{theorem}[\cite{3}, Theorem 4.1]
\begin{equation}\label{15.02.03-1}
\delta T = \beta_{a-}(x) + \beta_{a+}(x) + \beta_{d-}(x) + \beta_{d+}(x)\ge 0.
\end{equation}
Moreover, if $a$ is a strictly regular straight line strictly of rank one and if $a\cap d = \emptyset$, then $\delta T > 0$.
\end{theorem}

The proof given in \cite{3} for the case of $CAT(0)$-space remains for Busemann spaces without changes.

\subsection{Shadows}\label{Sha}

\begin{definition}[\cite{3}, Definition 4.3]
\emph{The complete shadow} of the point $x_0$ relatively the point $y \in \overline{X}\setminus \{x_0\}$ is by definition the set
\[\sha_y(x_0) := \{z \in \overline{X}|\ \exists [yz] \ x_0 \in [yz]\}.\]
Here the existence supposition is necessary only if both points $y,z$ are infinite: $y,z \in \partial_g X$.
\emph{The spherical shadow}  $\sha_y (x_0, \rho)$ of the point $x_0$ \emph{of radius} $\rho > 0$ relatively the point $y\in \overline{X}$ is the intersection of the shadow $\sha_y(x_0)$ with the sphere $S(x_0, \rho)$. In particular, if $\rho = +\infty$, then
\[\sha_y (x_0, +\infty) := \partial_g(\sha_y (x_0)) := \sha_y(x_0) \cap \partial_g X.\]
\end{definition}

\begin{theorem}\label{finite}
Let $x_0 \in X, y \in \overline{X}\setminus\{x_0\}$, and assume that if $y \in \partial_gX$ then $y$ is regular ideal point. Then for any numbers $\rho, \varepsilon > 0 $ there exists a number $\delta > 0$, such that if the point  $x_1 \in B(x_0, \delta)$ satisfies to equality $|yx_1| = |yx_0|$ (or $b_y(x_1) = b_y(x_0)$, for the case $y \in \partial_\infty X$), then
\[\sha_y(x_1, \rho) \subset \mathcal{N}_\varepsilon (\sha_y(x_0, \rho)).\]
\end{theorem}

\proof On the contrary, suppose that the claim is false: for some $\rho, \varepsilon > 0$ and any $\delta > 0$ there exist points $x_\delta \in B(x_0, \delta) \cap S(y, \, |yx_{0}|)$ and $z_\delta \in S(y, |yx_0|+\rho) \setminus \mathcal{N}_\varepsilon (\sha_y(x_0, \rho))$ such that $x_\delta \in [yz_\delta]$. 
In the case $y \in \partial_gX$, we have following changes: $x_\delta \in B(x_0, \delta) \cap \mathcal{HS}(\beta_y, x_{0})$ and $z_\delta \in \beta_{y}^{-1}(\rho) \setminus \mathcal{N}_\varepsilon (\sha_y(x_0, \rho))$, where $\beta_y$ is Busemann function corresponding to $y$ with $\beta_y(x_0) = 0$.
Fix a sequence $\delta_n \to 0$ and corresponding sequences of points $x_{\delta_{n}}$ and $z_{\delta_{n}}$. Obviously, $x_{\delta_{n}} \to x_0$. 
Since the space $X$ is proper, one can subtract converging subsequence from the sequence $z_{\delta_{n}}$. We assume that the sequence $z_{\delta_{n}}$ converges itself and 
$$\lim_{n \to \infty} z_{\delta_{n}} = z \in S(y, |yx_0|+\rho)$$
or
$$z \in \beta_y^{-1}(\rho)$$
when $y \in \partial_gX$.

For simplicity we finish the proof only for the case $y \in X$. The case $y \in \partial_gX$ is similar.
If $\gamma: |yz| \to X$ is natural parameterization of the segment $[yz]$, then
\begin{equation}\label{sered}
|x_0 \gamma(|yx_0|)| \le |x_0 x_{\delta_{n}}| + |x_{\delta_{n}} \gamma(|yx_0|)| \le \delta_n + |z_{\delta_{n}} z|
\end{equation}
The right hand in \eqref{sered} tends to zero when $n \to \infty$, hence the constant in the left hand is $0$. It follows that the point $z \in \sha_y(x_0, \rho)$ and points $z_n$ belong to $\mathcal{N}_\varepsilon (\sha_y(x_0, \rho))$ when $n$ is sufficiently large. 
This contradicts to their choice. \qed

The next statement simply follows from Theorem \ref{finite}.

\begin{corollary}\label{infty}

For any neighbourhood $\mathcal{N}_{y,\, K,\, \varepsilon}(\partial_g(\sha_y (x_0)))$  of the shadow at infinity $\partial_g (\sha_y(x_0))$ relatively $y \in X$ there exists a number $\delta>0$ such that for any point  $x_1\in B(x_0, \delta)$ the inclusion holds
\[\partial_g (\sha_y(x_1)) \subset \mathcal{N}_{y,\, K,\, \varepsilon} (\partial_g (\sha_y(x_0))).\]

\end{corollary}

Also we specify the situation for the case of strictly regular straight line strictly of rank one.

\begin{corollary}
Let $a: \R \to X$ be strictly regular straight line strictly of rank one. Suppose that both directions of the straight line $a$ in the point $x_0 = a(0)$  have unique opposite direction. Denote  $\xi = a(+\infty)$ and $\eta = a(-\infty)$. Then there exist numbers $K, \varepsilon > 0$ such that for any pair of points $\zeta \in \mathcal{N}_{x_{0},\, K,\, \varepsilon}(\sha_{\eta}(x_0))$ and $\theta \in \mathcal{N}_{x_0,\, K, \varepsilon}(\sha_{\xi}(x_0))$ the relation holds
$$\Td(\theta, \zeta) > \pi.$$
\end{corollary}

\proof Follows from the definition of the relation $\Td > 0$, Theorem \ref{finite} and Lemma \ref{12.04.08 8.04}. \qed

\subsection{Existence theorem for scissors}

In this section, we prove the following existence theorem.

\begin{theorem} \label{15.02.03-4}
Let $a: (-\infty , +\infty) \to X$ be strictly regular straight line strictly of rank one and  both directions of the straight line $a$ in the point $x_0 = a(0)$  have unique opposite direction. Then there exists a straight line $a'$ passing throw $a'(0)= x_0$ in the same directions with the following property. For any neighbourhood $\mathcal{U}$ of the triple
\begin{equation}\label{17.04.08 0.33}
(a'(+\infty), a'(-\infty), x_0) \in \partial_\infty X \times \partial_\infty X \times X
\end{equation}
there exist a triple $(\xi,  \eta, x) \in \mathcal{U}$ with $x \ne x_0$ and scissors $\langle a,b,c,d; x\rangle$ for which $b = [a(-\infty) \xi], c=[\eta a(+\infty)]$ and $d=[\eta\xi]$.
\end{theorem}

\begin{remark}\label{11.04.08 17.58}
Since the space $X$ is proper, any straight line $a$ in $X$ has infinite subset of points where both direction of $a$ have unique opposite direction (see \cite{3}, Theorem 4.3). Hence we can always choose appropriate parameterization for $a$.
\end{remark}

\proof 

Note that by the condition on the rank of the straight line $a$ every straight line $a'$ satisfying
\[a'(+\infty) \in \partial_g(\sha_{a(-\infty)}(x_0))\]
or
\[a'(-\infty) \in \partial_g(\sha_{a(+\infty)}(x_0))\]
has rank one. First we show that there exist scissors with base $a$ and the center arbitrarily close to $x_0$. By remark above, such consideration is also applicable to the straight line $a'$ passing throw $x_0$ in the direction of $a$, because in that case the lines $a$ and $a'$ are connected by asymptotic chain. 

Given a number $\rho > 0$ consider points $y' = a(-\rho)$ and $y'' = a(\rho)$. We have
\[\partial_g(\sha_{a(-\infty)}(x_0)) = \partial_g(\sha_{y'}(x_0))\]
and
\[\partial_g(\sha_{a(+\infty)}(x_0)) = \partial_g(\sha_{y''}(x_0)).\]

Furthermore, for some $K>\rho$ and $\varepsilon > 0$ there exist neighbourhoods
\[\mathcal{N}\,':=\mathcal{N}_{y', K, \varepsilon}(\partial_g(\sha_{y'}(x_0)))\]
and
\[\mathcal{N}\,'':=\mathcal{N}_{y'', K, \varepsilon}(\partial_g(\sha_{y''}(x_0)))\]
of shadows at infinity of the point $x_0$, such that
\begin{equation}\label{inftd}
Td(\xi, \eta) > \pi
\end{equation}
for each $\xi \in \mathcal{N}\,'$ and $\eta\in \mathcal{N}\,''$.

Choose $\delta_1$-neighbourhood $B(x_0, \delta_1)$ of the point $x_0$ defined from Corollary \ref{infty} within neighbourhoods $\mathcal{N}_{y', K, \varepsilon/2}(\partial(\sha_{y'}(x_0, \rho))$ and $\mathcal{N}_{y'', K, \varepsilon/2}(\partial(\sha_{y''}(x_0, \rho))$.

By Theorem \ref{finite} there exists $\delta_2$ such that for any point $x' \in B(x_0, \delta_2)$ inclusions hold:
\[\sha_{a(-\infty)}(x', \rho) \subset \mathcal{N}_{\varepsilon/2}(\sha_{a(-\infty)}(x_0, \rho))\]
and
\[\sha_{a(+\infty)}(x', \rho) \subset \mathcal{N}_{\varepsilon/2}(\sha_{a(+\infty)}(x_0, \rho)).\]

Denote $\delta_0 := \min \{\delta_1, \delta_2\}$. Then for any point $x \in \mathcal{U}_{\delta_{0}}(x_0)$ and straight lines $b'$ and $c'$ satisfying to conditions
\begin{itemize}
\item $b'(0) = c'(0) = x$,
\item $b'(-\infty) = a(-\infty)$ and
\item $c'(+\infty) = a(+\infty)$,
\end{itemize}
inclusions hold
\begin{equation}\label{bincl}
b'(\rho) \in \mathcal{N}_{\varepsilon/2}(\sha_{y'}(x_0, \rho))
\end{equation}
and
$$c'(-\rho) \in \mathcal{N}_{\varepsilon/2}(\sha_{y''}(x_0, \rho)).$$

We show that
\begin{equation}\label{b+infty}
b'(+\infty) \in \mathcal{N}\,'
\end{equation}
and
\begin{equation}\label{c-infty}
c'(-\infty) \in \mathcal{N}\,''.
\end{equation}

Let $\gamma: [0, +\infty) \to X$ be natural parameterization of the ray $\gamma = [y'b(+\infty)]$ and straight line $p$ passes throw $x_0 = (0)$ such that $p(+\infty) \in \partial_g(\sha_{y'}(x_0))$. Then
$$|\gamma(2\rho)p(\rho)| \le |\gamma(2\rho)b'(\rho)| +|b'(\rho)p(\rho)|.$$
The first item has an estimation
$$|\gamma(2\rho)b(\rho)|\le|\gamma(0)b'(-\rho)| = |a(-\rho)b'(-\rho)| \le |a(0)b'(0)| < \frac\varepsilon2.$$
By \eqref{bincl} the straight line $p$ can be chosen so that an estimation for the second item holds
$$|b(\rho)p(\rho)|<\frac\varepsilon2.$$
Finally
\[|\gamma(2\rho)p(\rho)|<\varepsilon,\]
proving the inclusion \eqref{b+infty}. The inclusion \eqref{c-infty} is analogous.

So we showed that there exists a straight line $d' \subset X$ connecting $c'(-\infty)$ and $b'(+\infty)$, forming scissors $\langle a,b',c',d'; x\rangle$.

Now, take a sequence $\delta_n \to 0$ and costruct for each $\delta_n$ scissors $\langle a_n, b_n, c_n, d_n; x_n \rangle$ with $|x_0x_n| < \delta_n$. Choose an accumulation triple point $(\xi',\eta', x_0) \in \partial_g X \times \partial_gX \times X$ for the sequence $(b_n(+\infty), c_n(-\infty), x_n)$. Then
$$\xi' \in \partial_g(\sha_{a(-\infty)}(x_0))$$ 
and 
$$\eta' \in \partial_g(\sha_{a(+\infty)}(x_0)).$$ 
Hence points $\xi'$ and $\eta'$ are connected by a straight line $a'=[\eta'\xi'] \subset X$ such that $a'(0)=x_0$. By the construction, the triple $(\xi',\eta', x_0)$ has needed neighbourhood $\mathcal U$. \qed

\subsection{Continuity of the shift function}

Here we prove that the shift function $\delta$ defined on appropriate subset in $\partial_\infty X \times \partial_\infty X \times X$ is continuous.

Let $a: \R \to X$ be strictly regular straight line strictly of rank one. Denote $Z(a)\subset \partial_\infty X \times \partial_\infty X \times X$ a subset consisting of triples $(\xi, \eta, x) \in \partial_\infty X \times \partial_\infty X \times X$ such that there exist scissors $\langle a,b,c,d; x\rangle$ with $b(+\infty) = \xi$ and $c(-\infty) = \eta$. 

For the completeness, we allow degenerate scissors. Scissors $\langle a, b, c, d; x \rangle$ are called {\it degenerate} if $x \in a \cap d$. We think that $(b(+\infty, c(-\infty), x) \in Z(a)$ as well.

\begin{theorem}\label{contin}
The shift function $\delta$ is continuous on the set $Z(a)$.
\end{theorem}

\proof
We use the equality \eqref{15.02.03-1}. Let the triple $(\xi_0, \eta_0, x_0) \in Z(a)$ be given. The point $x_0$ is the center of scissors $\langle a, b_0, c_0, d_0; x_0 \rangle$, where $b_0(+\infty) = \xi_0$ and $c_0(-\infty) = \eta_0$. 

Fix a number $\varepsilon > 0$. Then, by the continuity of Busemann functions $b_{a-}$ and $b_{a+}$ there exists a number $\sigma_1$, such that if the point $x'\in X$ satisfies inequality $|\, x_0x' | < \sigma_1$, then 
\begin{equation}\label{15.02.03-2}
|\, b_{a+}(x') + b_{a-}(x') - b_{a+}(x_0) - b_{a-}(x_0) | < 
\varepsilon/3.
\end{equation}

From the other hand, we use the regularity of all considering ideal points in $\partial_gX$. So, some neighbourhood $U$ of the point $\xi \in \partial_gX$ in geodesic ideal boundary is simultaneously the neighbourhood of horofunction ideal point represented by Busemann function $\beta_\xi$ and wise versa.

Fix a neighbourhood $V$ of the point $x_0$ with compact closure where values of Busemann functions $b_{d\pm}$ differ from $b_{d\pm}(x_0)$ at most by $\varepsilon/6$. Denote 
$$U_\pm (V) := \left\{f \in C(X)| \forall x \in V, |f - b_\pm(x_0)| < \frac16 \varepsilon \right\}$$
neighbourhoods of Busemann functions $b_\pm$ in $C(X)$ containing functions with values in $V$ differ from $b_\pm$ less then by $\varepsilon/6$. Let
$$\mathcal U_\pm := (U_\pm/\operatorname{consts}) \cap (\partial_mX) $$
be neighbourhoods in $\partial_mX = \partial_gX$ generated from $U_\pm (V)$. Regularity of points in $\partial_gX$ implies the following. If the ray $d'$ has endpoint $d'(+\infty) \in \mathcal U_+$ then Busemann function $\beta_{d'}$ behaves in the neighbourhood $V$ so that 
\begin{equation}\label{16.04.08 23.47}
|\, b_{d+}(x_0) - \beta_{d'}(x') - \const| < \varepsilon/3,
\end{equation}
for every $x'\in V$ and some constant. Similarly, if the ray $d''$ has endpoint $d''(+\infty) \in \mathcal U_-$, then
\begin{equation}\label{16.04.08 23.48}
|\, b_{d-}(x_0) - \beta_{d''}(x') - \const'| < \varepsilon/3
\end{equation}
for all $x' \in V$ and some constant $\const'$.
It follows from Lemma \ref{15.04.08 19.31} that for sufficiently small neighbourhoods $\mathcal U_+$ and $\mathcal U_-$ straight lines $d'$ and $d''$ passes arbitrarily close to each other, so the sum of constants $\const + \const'$ in \eqref{16.04.08 23.47} and \eqref{16.04.08 23.48} is arbitrarily close to $0$. In fact, all the considered neighbourhoods can be reduced so that the equality $\const = \const' = 0$ becomes admissible.

Denote
\[\mathcal U = (\mathcal U_+ \times \mathcal U_- \times V) \cap Z.\]
Then for any triple $(\xi', \eta', x') \in \mathcal U$ generating scissors $\langle a, b', c', d'; x'\rangle$ with the shift $\delta'$,
$$|\delta' - \delta| = |(b_{a-}(x_0) + b_{a+}(x_0) + b_{d-}(x_0) + b_{d+}(x_0)) -$$ $$-(b_{a-}(x') + b_{a+}(x') + b_{d'-}(x') + b_{d'+}(x'))| < \varepsilon.$$
Hence the claim. \qed

\section{Finish of the proof of main Theorem} 
Now we are ready to complete the proof of the Theorem \ref{th1.2}. 

Let $x, y \in X$ be two arbitrary points. There exists a straight line $a$ passing throw $x$ and $y$. As it was shown, $a$ is a straight line in the sense of both metrics $d_1$ and $d_2$. If $a$ is virtually of higher rank, then $d_1(x, y) = d_2(x, y)$ by Theorem \ref{17.04.08 0.25}. If $a$ is virtually singular, then $d_1(x, y) = d_2(x, y)$ by Corollary \ref{17.04.08 0.27}. Suppose now that $a$ is strictly regular and strictly has rank one. Then fix a point $x_0 = a(0)$ such that both directions of $a$ in $x_0$ have unique opposite. By Theorem \ref{15.02.03-4} there exists a straight line $a'$ passing throw $x_0$ in the same directions as $a$ such that for any neighbourhood $\mathcal{U}$ of the triple \eqref{17.04.08 0.33}
there exist a triple $(\xi,  \eta, x) \in \mathcal{U}$ with $x \ne x_0$ and scissors $\langle a,b,c,d; x\rangle$ for which $b = [a(-\infty) \xi]$, $c=[\eta a(+\infty)]$ and $d=[\eta\xi]$. Form degenerate scissors $\langle a, \bar b, \bar c, a', x_0 \rangle$, where straight line $\bar b$ is constructed from the ray $[a(-\infty), x_0]$ and the ray $[x_0 a'(+\infty)]$ and the straight line $c$ from the ray $[a'(-\infty) x_0]$ and the ray $[x_0 a(+\infty)$. The shift $\delta$ of degenerate scissors is equal to zero: $\delta = 0$. By Theorem \ref{contin} the shift $\delta(b(+\infty), c(-\infty), x)$ for scissors $\langle a, b, c, d, x \rangle$ closed to $\langle a, \bar b, \bar c, a', x_0 \rangle$ changes continuously and it tends to zero when $(b(+\infty), c(-\infty), x) \to (\bar b(+\infty), \bar c(-\infty), x_0)$.  Hence  its values cover some segment $[0, \Delta]$ where $\Delta > 0$.

Fix scissors $\langle a, b_n, c_n, d_n; x_n \rangle$ with
$$\delta(b_n(+\infty), c_n(-\infty), x_n) = \frac{1}{n}.$$
Then $n$-th degree of the transfer $T$ is an isometric translation on the distance $1$ along the straight line $a$ in the sense of both metrics $d_1$ and $d_2$. If $d_1(x, y) = k/n$, then $T^k(x) = y$ or $T^k(y) = x$. Since the images in the map $T$ does not depend on the choice of the metric $d_1$ or $d_2$, hence $d_2(x, y) = k/n$ as well. It follows that if $d_1(x, y)$ is rational, then $d_1(x, y) = d_2(x, y)$. Finally metrics $d_1$ and $d_2$ coincide because they are equivalent. \qed

\section{Some counterexamples}

\subsection{Trivial counterexamples}

Here we present some constructions leading to the counterexamples to the positive solution of A.D. Alexandrov problem. We begin with several elementary counterexamples.

\begin{counterexample}
Let the metric $d$ of the space $(X, d)$ does not takes values the unit distance. In particular, the diameter of the space $X$ can be less than $1$. Then every bijection of the space $X$ to itself presents the unit distance.
\end{counterexample}

\begin{counterexample}
The round Euclidean spheres $\S(o, \frac{1}{2\pi})$ and $\S(o, \frac{1}{\pi})$ in $\E^n$, $n \ge 2$. Both spheres admit the following bijection $\phi$. Let $A \subset \S$ be arbitrary proper centrally symmetric subset. Then we put
$$
\phi(x) = \left\{
	\begin{array}{ll}
		-x, & \mbox{ if } x \in A\\
		x & \mbox{otherwise}.
	\end{array}
\right.
$$
It is easy to see that in both cases the maps $\phi$ are bijections preserving the unit distance but not isometries.
\end{counterexample}

\begin{counterexample} The real line $\R$. The function 
$$f(x) = x + \frac{1}{2\pi}\sin (2\pi x)$$
preserves the unit distance, but it is not an isometry.
\end{counterexample}

\subsection{Grasshopper metric}

Let $(X, d)$ be arbitrary metric space. Here we construct new auxiliary metric $G_d$ on the set $X$ associated to the metric $d$.  We prove its property leading to several counterexamples for the positive answer to Alexandrov's problem in the class of $\R$-trees. The metric $G_d$ takes values in the set $\N \cup \{0, +\infty\}$.

\begin{definition}
Given points $x, y \in X$, {\it the grasshopper jump} of {\it length} $n \in \N$ from $x$ to $y$ is a map $j: \{0, \dots, n\} \to X$ such that $j(0) = x$, $j(n) = y$ and $d(j(i), j(i+1)) = 1$ for every $i \le n-1$. The grasshopper distance $G_d(x, y)$ is defined as minimal length of the grasshopper jump from $x$ to $y$. If there is no grasshopper jump from $x$ to $y$, we set $G_d(x, y) = +\infty$. The metric space $(X, d)$ is called {\it grasshopper jumps connected} if $G_d(x, y) < +\infty$ for all $x, y \in X$. {\it The grasshopper jumps component} of the point $x \in X$ is by definition the maximal grasshopper jump connected subspace of $X$ containing $x$. The subset $A \subset X$ is called {\it grasshopper invariant} if it is a union of grasshopper jumps components.
\end{definition}

\begin{lemma}\label{28.04.08 21.51}
Let the metric space $(X, d)$ be not grasshopper jumps connected and let $A \subset X$ be grasshopper invariant subset. Suppose that $A$ admits non-trivial isometry $\phi$ in the sense of the metric $G_d$ and there exists a point $y \in X \setminus A$ such that $d(y, z) \ne d(y, \phi(z)$ for some $z \in A$. Then the space $X$ admits a bijection $f: X \to X$ which preserves the distance 1 and is not an isometry.
\end{lemma}

\proof It is sufficient to take
$$
f(p) = \left\{
	\begin{array}{ll}
		\phi(p),& \mbox{ if } p \in A,\\
		p & \mbox{ otherwise.}
	\end{array}
\right.
$$
\qed

\begin{corollary}\label{28.04.08 21.54}
Let the space $X$ be a metric tree such that lengthes of all its edges are rational and presented by fractions with uniformly bounded denominators. Then $X$ admits non-isometric bijection preserving unit distance.
\end{corollary}

\proof Let $n$ be the largest common denominator of all fractions presenting lengthes of the edges in the tree $X$. We may assume that all the edges of $X$ have the same length $1/n$. Fix numbers $\alpha, \beta$ with $0 < \alpha < \beta < \frac{1}{2n}$. Let $A_\alpha \subset X$ (correspondingly, $A_\beta \subset X$) be the subset of all points in $X$ on the distance $\alpha$ (correspondingly $\beta)$ from the vertices. Then the sets $A_\alpha$, $A_\beta$ and 
$$A = A_\alpha \cup A_\beta$$
are grasshopper invariant and the metric spaces $(A_\alpha, G_d)$ and $(A_\beta, G_d)$ are isometric. The map $\phi_G: A \to A$ is defined by the following rule. Let the point $x \in A_\alpha$ lies in the edge $[a, b]$ and $|ax| = \alpha$. Then the image $y = \phi_G(x)\in A_\beta$ lies in $[a, b]$ and $|ay| = \beta$. Simultaneously we put $\phi_G(y) = x$. We show that the map $\phi_G$ is an isometry in the sense of the metric $G_d$ on $A$. First, note that $\phi_G$ is bijective, because the unique preimage of arbitrary point $y \in A_\beta$ is $x = A_\alpha \cap [ay]$ and the unique preimage of $x$ is $y$.

Let $p = x_0, x_1, \dots, x_n = q$ be a grasshopper jump from $p$ to $q \in A_\alpha$. Then $\phi(p) = y_0, y_1, \dots, y_n = \phi(q)$ is the grasshopper jump from $\phi(p)$ to $\phi(q)$. To see this it is sufficient to observe that all distances $|y_{i-1}y_i|$ for $i = \overline{1,n}$ are unit:
$$|y_{i-1}y_i| = 1.$$
Indeed, let $z_i, i= \overline{0,n}$ be the nearest vertex to $x_i$ and $y_i$. Then
$$|y_{i-1}y_i| = |x_{i-1}x_i| = |z_{i-1}z_i| = 1.$$
Hence
$$G_d(\phi(p), \phi(q)) \le G_d(p,q).$$
Analogously, 
$$G_d(p,q) \le G_d(\phi(p), \phi(q))$$
and hence $\phi$ preserves the grasshopper distance in $A$. Application the Lemma \ref{28.04.08 21.51} gives the result. \qed

\begin{remark}
Under the conditions of the Corollary \ref{28.04.08 21.54} it is easy to construct the continuous non-isometric bijection preserving the unit distance. In fact such bijection $\phi$ can be defined on edges by the formula
$$\phi(\gamma(t)) = \gamma\left(t + \frac{\sin 2\pi n t}{2 \pi n}\right),$$
where $\gamma: [0, \frac{1}{n}] \to X$ is the natural parameterization of the corresponding edge.
\end{remark}

\begin{corollary}
Let the $\R$-tree $X$ admits a non-trivial isometry and the set of distances between its branch points is at most countable. Then $X$ admits non-isometric bijection on itself preserving the unit distance.
\end{corollary}

\proof Let $\phi: X \to X$ be a non-trivial isometry and $P \subset \R_+$ be the set of distances between branch points of $X$. Fix a branch point $x \in X$. Let $\mathcal A$ be the grasshopper jumps component of the point $x$. Then $\mathcal A$ is a proper subset in $X$. Indeed, $\mathcal A$ is represented as the union
$$\mathcal A = \bigcup_{k=1}^\infty \mathcal A_k,$$
where
$$\mathcal A_k = \{y \in X|\ G_d(x, y) \le k\}.$$
Let the map $D_x: X \to \R_+$ be defined by the equality $D_x(y) = d(x, y)$. Then the set of values for the composition $D_x|_{\mathcal A}$ is at most countable, since  for each $k \in \N$ the set of values $D_x|_{\mathcal A_k}$ is at most countable. Hence $D_x|_{\mathcal A}$ does not contain any interval in $\R$. Hence the claim. Now we can apply the Lemma \ref{28.04.08 21.51}. \qed

\subsection{Maximum products}

\begin{definition}
Let $(X, d_X)$, $(Y, d_Y)$ be metric spaces. Then their {\it maximum product} is by definition the metric space $(X \times Y, d_{\max})$ where $d_{\max} = \max(d_X, d_Y)$ is the {\it maximum metric} on the product $X \times Y$:
$$d_{\max}((x_1, y_1), (x_2, y_2)) = \max\{d_X(x_1, x_2), d_Y(y_1, y_2)\}.$$ 
\end{definition}

\begin{counterexample}
Let the space $(Y, d_Y)$ admits non-isometric bijection $\phi$ preserving the unit distance. Then for any space $(X, d_X)$, the maximum product space $(X \times Y, d_{\max})$ also admits such bijection $\phi \times \Id$ acting by the formula 
$$(\phi \times \Id)(x, y) = (\phi(x), y).$$
The particular case is the maximum product $X \times \R$ for any space $(X, d_X)$.
\end{counterexample}

{\bf Acknowledgement}. The author is very grateful to S.V. Buyalo for the number of important corrections and remarks on the text. My special thanks to V.K. Kropina for very useful discussion on the subject of the paper.

\end{document}